\documentclass[11pt, a4paper]{article}

\usepackage{amssymb,amsmath}
\usepackage{cleveref}
\usepackage{graphicx}
\usepackage{url}
\usepackage{a4wide}
\usepackage{color}
\usepackage{algorithmic}
\usepackage{subcaption}

\definecolor{blue}{rgb}{0,0,.75}
\usepackage{listings} 

\lstnewenvironment{PseudoCode}[1][]
{\lstset{language=Matlab,basicstyle=\scriptsize, keywordstyle=\color{blue},numbers=left,xleftmargin=.04\textwidth,#1}}
{}

\usepackage{dutchcal}

\usepackage{multirow}

\title{A two-level  iterative scheme for general sparse linear systems based on approximate skew-symmetrizers\thanks{The first author was supported by Alexander von Humboldt Foundation for  a research stay at TU Berlin and the BAGEP Award of the Science Academy.
The second author was supported by Deutsche Forschungemeinschaft through collaborative research center SFB TRR 154 Project B03.
}}

\author{Murat Manguo\u{g}lu\footnotemark[2]
        \and Volker Mehrmann\footnotemark[3]}

\begin{document}

\maketitle

\renewcommand{\thefootnote}{\fnsymbol{footnote}}

\footnotetext[2]{Institut f\"{u}r Mathematik,  Technische Universit\"{a}t Berlin,  10623 Berlin,  Germany.
Present address: Department  of  Computer  Engineering,  Middle  East  Technical  University,  06800  Ankara, Turkey ({\tt manguoglu@ceng.metu.edu.tr}).}
\footnotetext[3]{Institut f\"{u}r Mathematik,  Technische Universit\"{a}t  Berlin,  10623 Berlin,  Germany ({\tt mehrmann@math.tu-berlin.de}).}

\begin{abstract}
We propose a two-level iterative scheme for solving general sparse linear systems. The proposed scheme consists of  a sparse preconditioner that  increases  the skew-symmetric part and makes the main diagonal of the coefficient matrix as close to the identity as possible. The preconditioed system is then solved via a particular Minimal Residual Method for Shifted Skew-Symmetric Systems ({\tt mrs}). This leads  to a two-level (inner and outer) iterative scheme where the {\tt mrs} has short term recurrences and satisfies an optimally condition. A preconditioner for the inner system is designed via a skew-symmetry preserving deflation strategy based on the skew-Lanczos process.  We demonstrate the robustness of the proposed scheme on sparse matrices from various applications.
\end{abstract}

\begin{keywords}
 symmetrizer, skew-symmetrizer, Krylov subspace method, shifted skew-symmetric system, skew-Lanczos method
\end{keywords}

\begin{AMS}
   	65F08, 65F10, 65F50
\end{AMS}

\section{Introduction}
We discuss the numerical solution of general linear systems
\begin{equation}
Ax=b,
\end{equation}
where $A\in\mathbb{R}^{n\times n}$ is a general large sparse invertible  matrix. If the coefficient matrix is symmetric and positive definite or symmetric and indefinite, one can use the Conjugate Gradient algorithm or the recently proposed two-level iterative scheme~\cite{manguoglu2019robust}, respectively. In this paper,  we propose a new  robust two-level black-box scheme for solving general systems  without any assumption on the symmetry or definiteness of the coefficient matrix.
In contrast to most other iterative methods, where preconditioning is often used to symmetrize the system and to lower the condition number, our new approach consists of an initial step which makes the system close to an identity plus skew-symmetric matrix that leads to  an effective shifted skew-symmetric preconditioner. Both the preconditioned system and the application of the preconditioner is approached by an iterative method  so that the method is a two-level (inner-outer) iterative scheme.

Our main motivation to study identity-plus-skew-symmetric preconditioners are linear systems arising in the time-discretization of dissipative Hamiltonian differential equations of the form
\begin{equation}\label{dissHS}
E\dot z= \left (J-R\right) z +f(t),\ z(t_0)=z_0
\end{equation}
where $\dot z$ denotes the derivative with respect to time, $J$ is a skew-symmetric matrix, $R$ is symmetric positive semi-definite  and $E$ is the symmetric positive semi-definite Hessian of a quadratic energy functional (Hamiltonian) $\mathcal H(z)=\frac 12 z^TEz$, see e.g. \cite{BeaMXZ18,Egg19,GraMQSW16,JacZ12,Sch13,SchJ14} for such systems in different physical domains and applications.  If one discretizes such systems in time, e.g. with the implicit Euler  method, and setting $z_k=z(t_k)$ then in each time step $t_k$ one has to solve a linear system of the form
\begin{equation}
(E-h(J-R))z_{k+1}= Ez_k+h f(t_k).
\end{equation}
Similar linear systems arise also when other discretization schemes are used.

The matrix $A=E+h(R-J)$ has a positive (semi)-definite symmetric part $M=E+hR$. If $M$ is positive definite, then with a two sided preconditioning with the Cholesky factor $L$  of $M=LL^T$, the matrix $L^{-1} A L^{-T}$ has the form $I+\tilde J$, where $\tilde J=hL^{-1} A L^{-T}$ is skew-symmetric {\color{red} \cite{concus1976generalized,Wid78}}. For such systems in \cite{concus1976generalized,IdeV07,Jia07,SzyW93,Wid78} structure exploiting Krylov subspace methods with three term recurrences were derived and analyzed.
Given a general square matrix $A$, symmetrizers from right or left, respectively, are symmetric  matrices, $S_r$ and $S_l$,  such that $AS_r=S^T_rA^T$ and $S_l A=A^TS^T_l$. Existing algorithms to construct dense and exact symmetrizers are studied and summarized in \cite{DopU16}. In this paper, however, we construct two-sided preconditioners  so that the preconditioned systems has the form $D+\hat J$, where $D$ is diagonal and close to the identity and $ \hat J$ is close to a skew-symmetric matrix (\emph{approximate shifted skew-symmetrizers (ASSS)}). To this preconditioned system we then apply a two-level iterative method, where the inner iteration is a skew-symmetric Krylov subspace method. We assume that both $A$ and $S\in\{S_r,S_l\}$ are sparse and nonsymmetric, with $S$ having a user defined sparsity structure.
The sparse ASSS preconditioner  is obtained by first applying a nonsymmetric permutation and scaling and then solving a sparse overdetermined linear least squares (LLS) problem  to obtain $S$.
Similar approaches for dense symmetrizers,  \cite{DopU16},
or algorithms for improving the {structural} symmetry in the context of sparse direct solvers, as proposed in  \cite{PorU16,Uce08}, do not have the latter property.

We note that while it is possible to obtain and use either  $S_r$  and $S_l$, in our experience the numerical results did not differ much as in left and right preconditioning. Therefore, in the rest of the paper we use the right variant and hereafter $S$ refers to $S_r$.

The paper is organized as follows. The proposed sparse approximate  skew-symmetrizer is introduced in \Cref{sec:alg},  a two-level Kyrlov subspace method based on the skew-symmetrizer is introduced in  \Cref{sec:prec}, numerical results are presented in \Cref{sec:results}, and  the conclusions follow in
\Cref{sec:conclusions}.

\section{A sparse approximate shifted skew-symmetrizing preconditioner}
\label{sec:alg}
Given a sparse invertible matrix $A\in \mathbb{R}^{n\times n}$, to achieve our goal of constructing a sparse approximate shifted skew-symmetrizing (ASSS) preconditioner,
we first apply  diagonal scalings  ($D_r,D_c$) and a row permutation  ($\mathcal{P}$),
\begin{equation}
\label{eq:mc64}
    \bar{A} = \mathcal{P}D_rAD_c
\end{equation}
such that the diagonal entries of $\bar{A}$ have modulus one and the off-diagonal elements are of modulus less than or equal to one. Such a permutation and scaling procedure is well established in the code  MC64 of the Harwell Subroutine Library (HSL)~\cite{hsl2007collection} and it is called the maximum product transversal with scaling. It solves a weighted bipartite matching problem and the resulting matrix  $\bar{A}$ is guaranteed to contain a zero-free main diagonal if $A$ is structurally nonsingular~\cite{DufK01}.
After this, we look for a sparse matrix $S$ such that
\begin{equation}
\label{eq1}
(\bar{A}S)_{i,j} = -(\bar{A}S)_{j,i}, \text{\quad for $i \neq j$},
\end{equation}
and
\begin{equation}
\label{eq2}
(\bar{A}S)_{i,i} = 1,  \text{\quad for $i=1,2,...,n$},
\end{equation}
where $S$ can have various sparsity structures, such as being  diagonal, tri-diagonal, banded, having the sparsity of $\bar{A}$, or any structure defined by the user.

The described problem can be formulated as a sparse over-determined LLS problem, where by \eqref{eq1}, each nonzero in the strictly upper triangular part of $|\bar{A}S|+|\bar{A}S|^T$ defines a constraint of the LLS problem and additional $n$ constraints  are obtained via \eqref{eq2}. Let $nz$ be the number of nonzeros in the strictly upper triangular part of $|\bar{A}S|+|\bar{A}S|^T$ and $nnz(S)$ be the number of nonzeros in $S$.  Then the LLS  problem has $nnz(S)$ unknowns and $nz+n$ equations, and if $nz+n > nnz(S)$ then the problem is overdetermined.

As a first example of a sparsity structure, let us assume $S=diag(s_{1,1},\ldots,s_{n,n})$, so that $nnz(S)=n$.  Then, \eqref{eq1} and \eqref{eq2} are given by
\begin{equation}
\label{eq:skew}
\bar{a}_{i,j}s_{j,j} + \bar{a}_{j,i}s_{i,i} = 0,
\end{equation}
and
\begin{equation}
\label{eq:diag}
\bar{a}_{i,i}s_{i,i} = 1,
\end{equation}
respectively. With $s =[s_{1,1}, s_{2,2}, ..., s_{n,n}]^T$, then
the resulting overdetermined system is given by
\begin{equation}
\label{lls_problem}
\begin{bmatrix}
B_u \\
B_l
\end{bmatrix}
s =
\begin{bmatrix}
\mathbf{0} \\
\mathbf{1}
\end{bmatrix}
\end{equation}
where $\mathbf{0}$ and $\mathbf{1}$ are vectors of all zeros of size $(nz+n)$ and all ones of size $n$, respectively. $B_u$ is a sparse matrix of size $(nz+n) \times n$, where each row has only two nonzeros,  $\bar{a}_{i,j}$ and $\bar{a}_{j,i}$ in its $i$-th and $j$-th columns, respectively, while $B_l$ is just the diagonal of $\bar{a}_{i,i}$. So with
\[
f(s):=
\left \|\begin{bmatrix}
B_u \\
B_l
\end{bmatrix}
s -
\begin{bmatrix}
\mathbf{0} \\
\mathbf{1}
\end{bmatrix} \right \|_2^2
\]
the unique solution of the  LLS problem is obtained by computing $\min_{s}{f(s)}$.
The unique solution can be obtained via a direct or iterative sparse LLS  solver. In order to obtain more flexibility in the importance of the two constraints, we introduce a weighting parameter ($\gamma>0$), and solve the weighted problem
\begin{equation}
\label{eq:regularized}
f(s,\gamma) =
\left \|  B_u s  \right \|_2^2  + \gamma \left \| B_l s - \mathbf{1}\right \|_2^2= \left \|(\bar{A}S) + (\bar{A}S)^T \right \|_F^2  + \gamma \left \|\mathcal{D}(\bar{A}S) - I \right \|_F^2,
\end{equation}
where $\mathcal{D}(X)$ denotes a diagonal matrix whose diagonal entries are those of $X$.

For a general sparse $S$, the LLS problem is formulated in a similar way as in the diagonal case. The set of constraints is defined for each nonzero $(i,j)$ in the strictly upper (or lower) triangular nonzero pattern of the matrix $|\bar{A}S|+|\bar{A}S|^T$ via \eqref{eq1}, (using Matlab column notation) via
\begin{equation}
\bar{A}_{i,:}S_{:,j} + \bar{A}_{j,:}S_{:,i} = 0,
\end{equation}
and the diagonal constraints are obtained for $i=1,2,...,n$ via \eqref{eq2}.
Note that one needs to map the nonzero entries of $S$ to a vector in order to form the LLS problem and map it back to $S$ after obtaining the solution vector.  This can be done using the sparse matrix storage format.
In Appendix~\ref{app:A}, we present a Matlab implementation which stores the non-zeros of sparse matrices in column major order, i.e. compressed sparse column format.

\section{ A bilevel iterative scheme}
\label{sec:prec}
Given a {general} sparse linear system
\begin{equation}
\label{eq:sys}
Ax=b
\end{equation}
where $A\in \mathbb{R}^{n\times n}$ is nonsingular. As discussed in the introduction, many preconditioners are either applied or aim for a  symmetric or symmetric positive definite system, since for these we have short recurrences in Krylov subsapce methods like the conjugate gradient method. Only very few algorithms focus on  skew-symmetric or shifted skew-symmetric  structure. In this section we present the theoretical basis for an algorithm that preprocesses the system such that the coefficient matrix is as close as possible to a shifted skew-symmetric matrix and then use the shifted skew-symmetric part of the matrix as preconditioner applying it as an iterative solver with short recurrences and optimality property that requires only one inner product per iteration.

Consider   the splitting of the coefficient matrix into its symmetric and skew-symmetric part
\begin{equation}
A= M+J
\end{equation}
where $M=M^T$ and $J=-J^T$, (in applications often  $J$ is even a matrix of small norm). If $M$ is positive definite then one can precondition the system by computing the  Cholesky factorization $M=LL^T$ and solve the modified system
\begin{equation}
(I+L^{-1}JL^{-T})L^Tx = L^{-1}b
\end{equation}
where $L^{-1}JL^{-T}$ is again skew-symmetric.
However, in general $M$ is not positive definite, it may be indefinite or singular. In this case we propose a black-box algorithm that employs an ASSS preconditioner. This is a two-level procedure, in which we first apply the discussed nonsymmetric row permutation and scaling to obtain a zero free diagonal with diagonal entries of modulus one and  off-diagonal entries of modulus less than or equal to one. The second step applies a sparse matrix $S$ obtained via the algorithm described in Section \ref{sec:alg} by solving a sparse LLS problem. After ASSS preconditioning, the modified system is has the form
\begin{equation}
\label{eq:modified}
\widehat{A} \widehat{x}  = \widehat{b}
\end{equation}
where $\widehat{A} = \mathcal{P}D_rAD_cS$, $\widehat{x} = S^{-1}D_c^{-1}x $ and $\widehat{b}  =  \mathcal{P}D_rb$.  Let
$
\widehat{M} = \frac{\widehat{A}+\widehat{A}^T}{2}$, and $\widehat{J} = \frac{\widehat{A}-\widehat{A}^T}{2}$
Note that due to the ASSS preconditioning, even though  $\widehat{M}$ is still not guaranteed to be positive definite, it has eigenvalues clustered around $1$  and typically very few negative eigenvalues. Furthermore, the skew-symmetric part, $\widehat{J}$,  is more dominant now. One can now compute a Bunch-Kaufman-Parlett factorization \cite{BunKP76}, $\widehat{M} = \widehat{L}\widehat{D}\widehat{L}^T$ and  modify the factorization to obtain
\begin{equation}\label{absfac}
|\widehat{M}| = \widehat{L}|\widehat{D}|\widehat{L}^T
\end{equation}
where,  as in \cite{vecharynski2013absolute},  $|\widehat{D}|=V|\Lambda|V^T$ if $\widehat{D}$ has a spectral decomposition $V\Lambda V^T$. Then, $|\widehat{D}|$ has a Cholesky factorization  $L_{|\widehat{D}|}L_{|\widehat{D}|}^T$, since it is positive definite. Setting $\mathcal{L} := \widehat{L}L_{|\widehat{D}|}$, and multiplying  \eqref{eq:modified} from the left with $\mathcal{L}^{-1}$ and inserting $I=\mathcal{L}^{-T}\mathcal{L}^{T}$, we obtain the system
$ \mathcal{A} \mathcal{x} = \mathcal{b} $, where  $\mathcal{A} = \mathcal{L}^{-1}\widehat{A} \mathcal{L}^{-T}$, $\mathcal{x} = \mathcal{L}^T \widehat{x}$ and $\mathcal{b} =\mathcal{L}^{-1} \widehat{b}$.  We note that $\mathcal{A}$ can be split as
\begin{equation}
\mathcal{A} = \underbrace{(L^{-1}_{|\widehat{D}|}\widehat{D} L^{-T}_{|\widehat{D}|} -I )}_{\mathcal{M}_r}  + (I + \underbrace{\mathcal{L}^{-1}\widehat{J}\mathcal{L}^{-T}}_{\mathcal{J}})
\end{equation}
where the rank of $\mathcal{M}_r$ is   equal to the number of negative eigenvalues of $\widehat{D}$ which is expected to be very small  and $I+\mathcal{J}$ is a shifted skew-symmetric matrix. Furthermore, $\mathcal{M}_r$ is symmetric and block diagonal with only a few nonzero blocks of size either $1\times 1 $ or $2 \times 2$ and is of rank $r \ll n$. The $1\times 1$ blocks have the value $-2$ and the $2\times 2$ blocks have eigenvalues $\{-2,0\}$.  Due to the (almost) diagonal and low rank structure of $\mathcal{M}_r$, it is simple to obtain a symmetric low-rank decomposition
\begin{equation}
\label{eq:low_rank_decomp}
\mathcal{M}_r = U_r \Sigma_r U_r^T,
\end{equation}
where $\Sigma_r=-2 I_r$, and $U_r$ is a sparse (with either one or two nonzero entries  per column) $n\times r$ orthogonal matrix. A pseudocode for computing such low rank decomposition is presented in Figure~\ref{fig:low_rank_decomposition}.

\begin{figure}[ht]
\begin{algorithmic}
{
\STATE {\bf Input:} $\mathcal{M_r}\in \mathbb{R}^{n\times n}$, $r$ (rank of $\mathcal{M}_r)$,  $\mbox{\rm Ind}$ (set of indices of nonzeros of $\mathcal{M}_r$)
\STATE $\Sigma_r \gets 0,U_r \gets 0,U \gets 0$
\STATE $i \gets 1, j \gets 1$
\STATE $\mathcal{M}'_r \gets \mathcal{M_r}(\mbox{\rm Ind},\mbox{\rm Ind})$
\WHILE{($i<r$)}
\IF{($\mathcal{M}'_r(i,i+1) =0 )$}
\STATE $\Sigma_r(j,j) \gets \mathcal{M}'_r(i,i)$
\STATE $U(:,j)  \gets e_i$
\STATE $i \gets i+1$
\ELSE
\STATE {Compute the eigenpair \{$\lambda_2,v_2$\} of  $\mathcal{M}'_r(i:i+1,i:i+1)$}
\STATE $\Sigma_r(j,j) \gets \lambda_2$
\STATE $U(:,j) \gets [e_i, e_{i+1}]v_2$
\STATE $i \gets i+2$
\ENDIF
\STATE $j \gets j +1$
\ENDWHILE
\IF{(i=r)}
\STATE $\Sigma_r(j,j) \gets \mathcal{M}'_r(i,i)$
\STATE $U(:,j) \gets e_i$
\ENDIF
\STATE $U_r(\mbox{\rm Ind},:)\gets U$
\STATE {\bf Output:} $U_r \in \mathbb{R}^{n\times r}, \Sigma_r \in \mathbb{R}^{r \times r}$
}
\end{algorithmic}
\caption{Sparse low rank decomposition of $\mathcal{M}_r=U_r\Sigma_rU_r^T$.}
\label{fig:low_rank_decomposition}
\end{figure}

The cost of this last step is $O(r)$ arithmetic operations, since it only needs to work with a submatrix of $\mathcal{M}_r$ corresponding to indices of nonzero entries. Using this factorization, we obtain
$\mathcal{A} = U_r \Sigma_r U_r^T + {\mathcal{S}}$
where $\mathcal{S}  = I + \mathcal{J}$, so that $\mathcal{A}$ is a shifted skew-symmetric matrix with a low-rank perturbation.
Using the Sherman-Morrison-Woodbury formula \cite{GolV96}, we theoretically have the exact inverse
\begin{equation}
    \mathcal{A}^{-1} = \mathcal{S}^{-1} - \mathcal{S}^{-1}U_r (\Sigma^{-1}_r + U_r^T \mathcal{S}^{-1}U_r)^{-1}U_r^T \mathcal{S}^{-1},
\end{equation}
which can be applied to the right hand side vector $\mathcal{b}$ to obtain $\mathcal{x}$, by solving only shifted skew-symmetric linear systems.

In practice, for large scale sparse systems, it is typically too expensive and storage intensive to compute the full  $LDL^T$ factorization of $\widehat{M}$, instead, an incomplete factorization $\tilde{M}  = \tilde{L}\tilde{D}\tilde{L}^T$ can be utilized together with the Cholesky factorization of $|\tilde{D}| = L_{|\tilde{D}|}L_{|\tilde{D}|}^T$ where $\mathcal{\tilde{L}} = \tilde{L}L_{|\tilde{D}|}$. This leads to a modified system,
\begin{equation}
\tilde{\mathcal{L}}^{-1}\widehat{A}\tilde{\mathcal{L}}^{-T}\tilde{\mathcal{L}}^T \widehat{x} = \tilde{\mathcal{L}}^{-1}\widehat{b}
\end{equation}
which then is solved iteratively using a Krylov subspace method with the preconditioner
\begin{equation}
\label{preconditioner}
P = \underbrace{(L^{-1}_{|\tilde{D}|}\tilde{D} L^{-T}_{|\tilde{D}|} -I )}_{\mathcal{\tilde{M}}_r}  + (I + \underbrace{\mathcal{\tilde{L}}^{-1}\widehat{J}\mathcal{\tilde{L}}^{-T}}_{\mathcal{\tilde{J}}})
\end{equation}
or alternatively, if $r$ is  zero, with a preconditioner
\begin{equation}
\label{preconditioner2}
P = \tilde {\mathcal S}=I + \mathcal{\tilde{J}}.
\end{equation}

One can, in principle, even apply $P^{-1}$ exactly as  described earlier. However, in a practical implementation applying $\tilde{\mathcal{S}}^{-1}$ via a direct solver is expensive, therefore one can apply it approximately by solving  shifted skew-symmetric systems iteratively, where the coefficient matrix is $\tilde{\mathcal{S}}$. This gives rise to an inner-outer iterative scheme.
In addition to solving a  shifted skew-symmetric system (where the coefficient matrix does not have to be formed explicitly) with a single right hand side vector, applying $P^{-1}$ requires sparse matrix-vector/vector-vector operations, solution of a dense $r\times r$ system, as well as one time cost of computing a low rank decomposition of $\mathcal{\tilde{M}}_r=\tilde{U}_r \tilde{\Sigma}_r \tilde{U}_r^T $ and  solving a shifted skew-symmetric system with multiple right hand side vectors.
The convergence rate of the outer Krylov subspace method depends on the spectrum of the preconditioned coefficient matrix $P^{-1}\tilde{\mathcal{L}}^{-1}\widehat{A}\tilde{\mathcal{L}}^{-T}$. The incomplete factorization of $\widehat{M}$ is an approximation such that $\widehat{M}=\tilde{L}\tilde{D}\tilde{L}^T + E$, where $E$ is a small norm error matrix. Assuming we apply $P^{-1}$ exactly, then the preconditioned coefficient matrix is  $P^{-1}\tilde{\mathcal{L}}^{-1}\widehat{A}\tilde{\mathcal{L}}^{-T}=I+P^{-1}\tilde{\mathcal{L}}^{-1}E\tilde{\mathcal{L}}^{-T}$. Due to the sparse ASSS preconditioning step, $\widehat{M}$ is already  close to identity and $\widehat{J}$ is dominant. Therefore, the norm of the perturbation of the preconditioned matrix from identity  ($||P^{-1}\tilde{\mathcal{L}}^{-1}E\tilde{\mathcal{L}}^{-T}||$) is expected to be small.

\subsection{Solution of  sparse shifted skew-symmetric sytems}
Application of the described  preconditioners involve the solution of linear systems where the coefficient matrix, $I + \mathcal{\tilde{J}}$, is shifted skew-symmetric. Specifically, we are interested in the iterative solution of such systems. While general algorithms such as Bi-Conjugate Gradient Stabilized ({\tt bicgstab}) \cite{Van92}, Generalized Minimal Residual ({\tt gmres}) \cite{SaaS86}, Quasi Minimal Residual ({\tt qmr}) \cite{FreN91} and Transpose Free Quasi Minimal Residual ({\tt tfqmr}) \cite{Fre93} can be used,  there are some iterative solvers available for shifted skew-symmetric systems such as   {\tt CGW} \cite{concus1976generalized,SzyW93,Wid78} and the Minimal Residual Method for Shifted Skew-Symmetric Systems ({\tt mrs}) \cite{IdeV07,Jia07}. We use {\tt mrs}, since it has a short recurrence and satisfies an optimality property. Furthermore, {\tt mrs} requires only one inner product per iteration~\cite{IdeV07} which would be  a great advantage if the algorithm is implemented in parallel since inner products require all to all reduction operations which create a synchronization points.  In addition to shifted skew-symmetric systems with one right-hand side vector, we also need to solve such systems with multiple right-hand side vectors. As far as we know, currently, there is no  "block" {\tt mrs} available.
Even though block Krylov methods are more amenable to breakdown, there are also ways to avoid the break down (for example, for block-CG see \cite{ji2017breakdown}). We  instead implemented a version of the {\tt mrs}  algorithm based on simultaneous iterations for multiple right-hand side vectors which is given in Appendix~\ref{app:B}.
In the proposed scheme, the convergence rate of {\tt mrs} iterations depends on the spectrum of the shifted skew-symmetric coefficient matrix, $I+\tilde{\mathcal{J}}$. In the next subsection, we propose a technique to improve this spectrum while preserving its shifted skew-symmetry.

\subsection{Improving the spectrum of shifted skew-symmetric systems via deflation}
One disadvantage of the {\tt mrs} algorithm is that if a preconditioner is used, then the preconditioned system should be also shifted skew-symmetric which may not be easy to obtain. Therefore, we propose an alternative deflation strategy to improve the number of iterations of {mrs}. For a shifted skew-symmetric system,
\begin{equation}
\label{eq:ss}
(I+\mathcal{\tilde{J}})z = y.
\end{equation}
we eliminate the extreme eigenvalues of $I+\mathcal{\tilde{J}}$, by running $k$-iterations ($k\ll n) $ of the skew-Lanczos process on $\mathcal{\tilde{J}}$, see \cite{GreV09,IdeV07,Jia07}. A pseudocode for this procedure is presented in Figure~\ref{fig:skewlanc}.
\begin{figure}[ht]
\begin{algorithmic}
{
\STATE {\bf Input:} $\mathcal{\tilde{J}}\in \mathbb{R}^{n\times n}$ ($\mathcal{\tilde{J}} = - \mathcal{\tilde{J}}^T$) and $k$.
\STATE Let $q_1$ be an arbitrary vector $\in \mathbb{R}^n$
\STATE $q_1 \gets q_1/||q_1||_2$
\STATE $z \gets \mathcal{\tilde{J}} q_1$
\STATE $\alpha_1 \gets ||z||_2$
\IF {$\alpha_1 \neq 0$}
    \STATE $q_2 \gets -z/\alpha_1$
    \FOR{$i=2$ to $k-1$}
        \STATE $z \gets \mathcal{\tilde{J}}q_{i}-\alpha_{i-1}q_{i-1}$
        \STATE $\alpha_i \gets ||z||_2$
        \IF {$\alpha_i = 0$}
            \STATE break
        \ENDIF
        \STATE $q_{i+1} \gets -z/\alpha_i$
    \ENDFOR
\ENDIF
\STATE {\bf Output:} $Q_k = [q_1, q_2, ... , q_k]\in \mathbb{R}^{n\times k}$ and $\tau = [\alpha_1,\alpha_2,...,\alpha_{k-1}]^T \in \mathbb{R}^{k-1}$
}
\end{algorithmic}
\caption{Skew-Lanczos procedure}
\label{fig:skewlanc}
\end{figure}
Considering the resulting matrices
\begin{equation}
S_k =
\begin{bmatrix}
0 & \alpha_1 & & 0\\
-\alpha_1 & \ddots & \ddots &  \\
& \ddots & \ddots &  \alpha_{k-1}
\\0 &  & -\alpha_{k-1} & 0
\end{bmatrix}, \quad
Q_k =
\begin{bmatrix}
q_1, q_2, \hdots , q_k,
\end{bmatrix}
\end{equation}
we deflate the system in \eqref{eq:ss} by forming
\begin{equation}
[(I+\mathcal{\tilde{J}})-Q_kS_kQ_k^T + Q_kS_k Q_k^T ]z = y,
\end{equation}
such that $Q_k^T\mathcal{\tilde{J}}Q_k=S_k$ where $Q_k$ is $n\times k$ with $Q_k^TQ_k = I$ and $S_k$ is a tridiagonal skew-symmetric $k\times k$ matrix. Let $\bar{\mathcal{J}}=\mathcal{\tilde{J}}-Q_kS_kQ_k^T$ which is still skew-symmetric and the largest (in modulus) eigenvalues have been set to zero. Then the system in \eqref{eq:ss} can be written as a low rank perturbation of a  shifted skew-symmetric system
\begin{equation}
\label{eq:deflated_shifted_skewsymmetric_system}
[(I+\bar{\mathcal{J}}) + Q_kS_kQ_k^T]z = y
\end{equation}
which can be handled again by the Sherman-Morrison-Woodbury formula. In fact, this low rank perturbation can be combined with the low rank perturbation in \eqref{preconditioner}, i.e., the preconditioner $P$ can be rewritten as
\begin{equation}
\label{eq:preconditioner}
P =\begin{bmatrix}
Q_k, \tilde{U}_r \end{bmatrix} \begin{bmatrix}
S_k &   \\
  & \tilde{\Sigma}_r
\end{bmatrix}
\begin{bmatrix}
Q^T_k \\
 \tilde{U}_r^T
\end{bmatrix}
+\bar{\mathcal{S}}
\end{equation}
where $\bar{\mathcal{S}} = I + \bar{\mathcal{J}}$. Then,  the preconditioner  can  be applied directly as via
\begin{equation}
\label{eq:smw}
P^{-1} = \bar{\mathcal{S}}^{-1} - \bar{\mathcal{S}}^{-1} \bar{U}_{r+k} (\bar{\Sigma}_{r+k} + \bar{U}^{T}_{r+k}\bar{\mathcal{S}}^{-1}\bar{U}_{r+k})^{-1}\bar{U}^{T}_{r+k}\bar{\mathcal{S}}^{-1}
\end{equation}
where $\bar{U}_{r+k}=\begin{bmatrix}
Q_k, \tilde{U}_r \end{bmatrix}$ and $\bar{\Sigma}_{r+k}=\begin{bmatrix}
S_k &   \\
  & \tilde{\Sigma}_r
\end{bmatrix}$.
Note that $P$ is the same preconditioner as in \eqref{preconditioner}, except that the perturbation is of rank $r+k$ now and the shifted skew-symmetric matrix ($\bar{\mathcal{S}}$) has a better spectrum, see Section~\ref{sec:deflation}.

\section{Numerical results}
\label{sec:results}

\subsection{Implementation details for the numerical experiments} As a baseline of comparison, we implemented a robust general iterative scheme that was proposed in \cite{BenHT00}.  It uses the same permutation and scalings given in (\ref{eq:mc64}) followed by a symmetric permutation.  We use Reverse Cutthill-McKee (RCM) reordering since RCM reordered matrices have  better robustness in subsequent applications of ILU type preconditioners \cite{BenHT00}.  After the symmetric permutation, we use ILU preconditioners with no fill-in ($ilu(0)$),  with pivoting and threshold of $10^{-1}$ ($ilutp(10^{-1})$) and $10^{-2}$ ($ilutp(10^{-2})$) of Matlab.  We call this method  as {\tt mps-rcm} and it is implemented in Matlab R2018a.

Our new method is also implemented in Matlab R2018a in two stages: preprocessing  and iterative solution.  In the preprocessing stage, we obtain  the sparse ASSS preconditioner where we just need the coefficient matrix to obtain  the permutation and scalings by calling HSL-MC64 via its Matlab interface. Followed by solving the LLS problem in (\ref{lls_problem}), which we do directly via Matlab's backslash operation. Then, we compute an incomplete Bunch-Kaufman-Parlett factorization of $\widehat{M}$ via the Matlab interface of  $sym\text{-}ildl$ software package~\cite{greif2017sym}. We use its default parameters except we disable any further scalings.  Similar to $ilutp$, we use  $10^{-1}$ and $10^{-2}$ thresholds and allow any fill-in and similar to $ilu(0)$, we  allow as many nonzeros as the original matrix per column with no threshold based dropping. We call these:  $ildl(10^{-1})$, $ildl(10^{-2})$ and $\sim ildl(0)$, respectively. We compute the low rank factorization in (\ref{eq:low_rank_decomp}) and apply a few steps of the skew-Lanczos process to deflate the shifted skew-symmetric part of the coefficient matrix. Finally, we iteratively solve the shifted skew-symmetric linear system of equations that arise in (\ref{eq:smw}) with multiple right hand side vectors via {\tt mrs}
\begin{equation}
    \bar{\mathcal{S}}X = \bar{U}_{r+k}
\end{equation}
and form the $(r+k)\times (r+k)$ dense matrix
\begin{equation}
\bar{\Sigma}_{r+k} + \bar{U}^{T}_{r+k}\bar{\mathcal{S}}^{-1}\bar{U}_{r+k}
\end{equation}
explicitly. All of these preprocessing steps do not require the right hand side vector and they are done only once if a sequence of linear systems with the same coefficient matrix but with different right hand side vectors need to be solved.

After preprocessing, the linear system of equation in (\ref{eq:modified}) is solved via a Krylov subspace method with the preconditioner in (\ref{eq:preconditioner}). At each iteration of the Krylov subspace method, the inverse of the preconditioner is applied as in (\ref{eq:smw}). This requires the solution of a shifted skew-symmetric linear system. We use the {\tt mrs} method for those shifted skew-symmetric systems.

 As the outer Krylov subspace method, some alternatives are {\tt bicgstab}, {\tt gmres}, and {\tt tfqmr}. Even though they often behave almost the same~\cite{BenHT00}, {\tt gmres} requires a restart parameter that defies our objective toward obtaining  a {\em black-box} solver and {\tt bicgstab} has erratic convergence. Alternatively,  {\tt tfqmr} has a smoother convergence and does not require restarting. We observe that {\tt tfqmr} can stagnate, which is also noted in~\cite{zhang2000preconditioned}. Therefore, as a challenge for our new approach, we use {\tt tfqmr} for both our proposed scheme and {\tt mps-rcm}. The stopping criterion for {\tt tfqmr} is set to $10^{-5}$ and for the inner {\tt mrs} iterations of the proposed scheme, we use the same stopping criterion. The right hand side is determined from the solution vector of all ones.

 We note that even though we use Matlab's built-in functions as much as possible while implementing the proposed iterative scheme,  {\tt mps-rcm} is entirely using the built-in functions of Matlab or efficient external libraries. Therefore, no fair comparison in terms of the running times in Matlab is currently possible. An efficient and parallel implementation of the proposed scheme requires a lower level programming language such as C/C++ due to the low-level algorithmic and data structural details that need to be addressed efficiently. For example, the proposed scheme needs efficiently accessing rows and columns of a sparse matrix. At first glance, one might tend to store the matrix both in Compressed Sparse Row and Column formats, however this approach is doubling the memory requirements. Therefore, a new storage scheme without much increase in the memory requirements is needed. Also, efficient and parallel implementation of sparse matrix-vector multiplications, where the coefficient matrix is symmetric (and shifted skew-symmetric) and parallel sparse triangular backward/forward sweeps are  challenging problems. These are still active research areas by themselves~\cite{alappat2020recursive,ccuugu2020parallel}. Therefore, we leave these issues as future work and focus on the robustness of the proposed scheme in Matlab.

\subsection{Test problems}
\begin{figure}[ht]
\begin{subfigure}{.32\textwidth}
  \centering
  \includegraphics[width=1\linewidth]{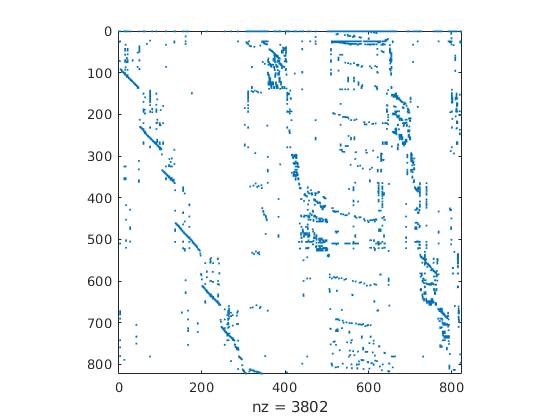}
  \caption{bp\_200}
  \label{fig:spy:bp_200}
\end{subfigure}%
\begin{subfigure}{.32\textwidth}
  \centering
  \includegraphics[width=1\linewidth]{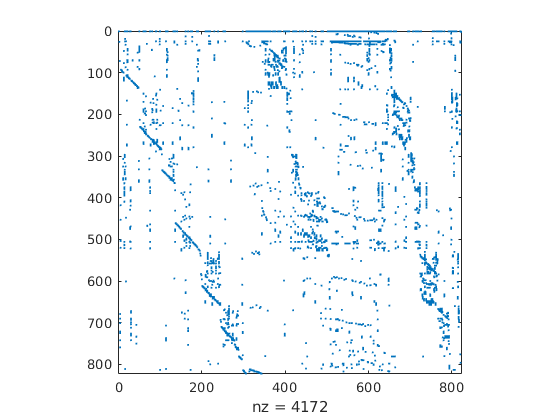}
  \caption{bp\_600}
  \label{fig:spy:bp_600}
\end{subfigure}
\begin{subfigure}{.32\textwidth}
  \centering
  \includegraphics[width=1\linewidth]{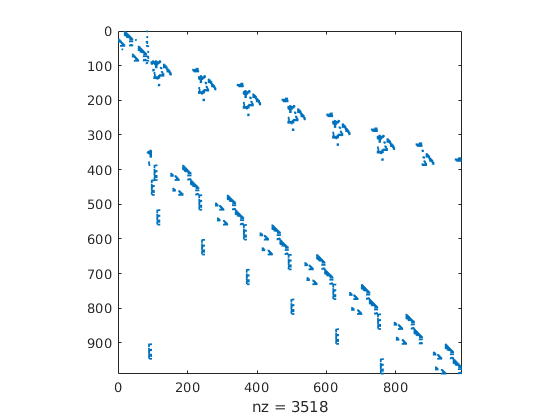}
  \caption{west0989}
  \label{fig:spy:west0989}
\end{subfigure}
\begin{subfigure}{.32\textwidth}
  \centering
  \includegraphics[width=1\linewidth]{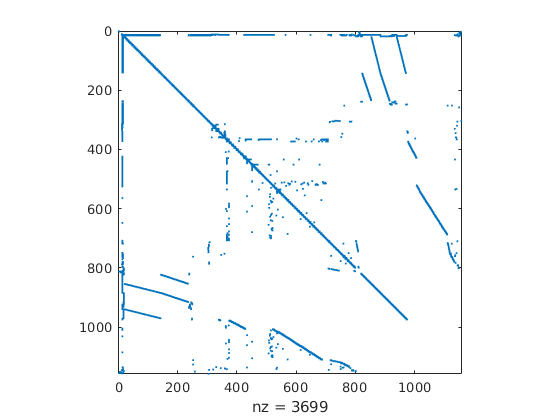}
  \caption{rajat19}
  \label{fig:spy:rajat19}
\end{subfigure}
\begin{subfigure}{.32\textwidth}
  \centering
  \includegraphics[width=1\linewidth]{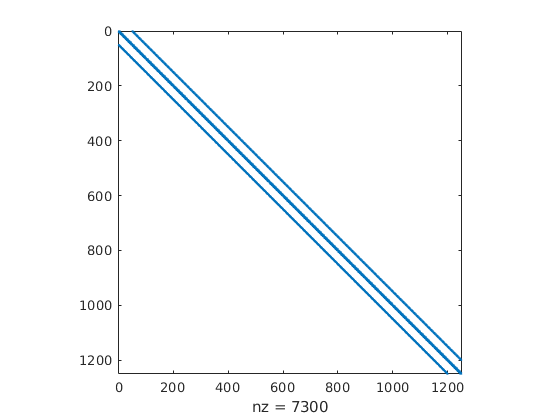}
  \caption{rdb1250l}
  \label{fig:spy:rdb1250l}
\end{subfigure}
\begin{subfigure}{.32\textwidth}
  \centering
  \includegraphics[width=1\linewidth]{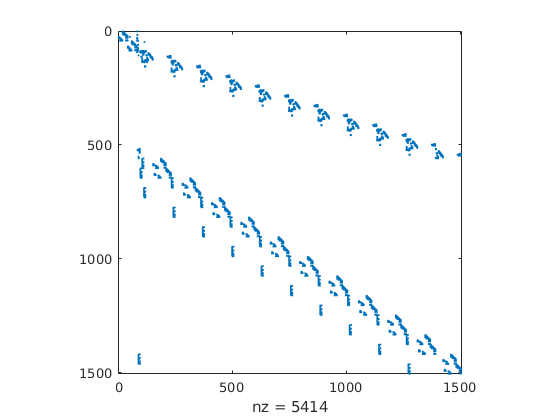}
  \caption{west1505}
  \label{fig:spy:west1505}
\end{subfigure}
\begin{subfigure}{.32\textwidth}
  \centering
  \includegraphics[width=1\linewidth]{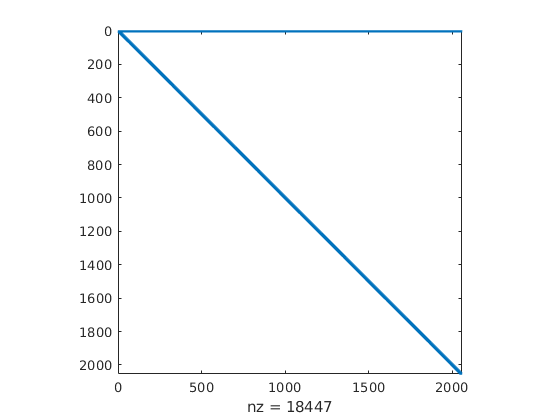}
  \caption{chebyshev2}
  \label{fig:spy:chebyshev2}
\end{subfigure}
\begin{subfigure}{.32\textwidth}
  \centering
  \includegraphics[width=1\linewidth]{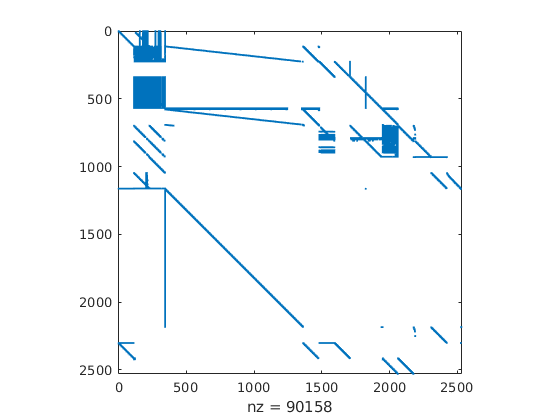}
  \caption{orani678}
  \label{fig:spy:orani678}
\end{subfigure}
\begin{subfigure}{.32\textwidth}
  \centering
  \includegraphics[width=1\linewidth]{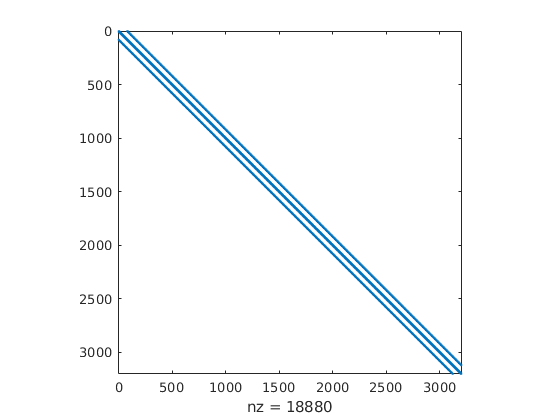}
  \caption{rdb3200l}
  \label{fig:spy:rdb3200l}
\end{subfigure}
\begin{subfigure}{\textwidth}
  \centering
  \includegraphics[width=0.32\linewidth]{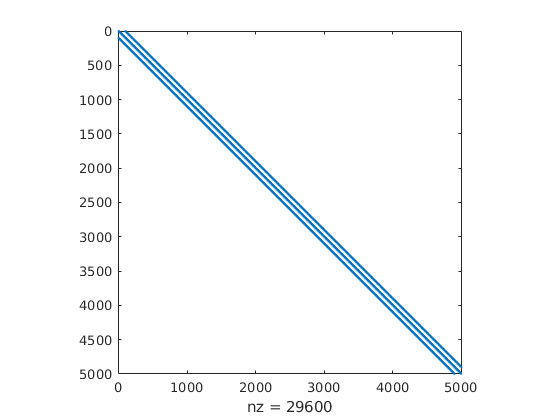}
  \caption{rdb5000}
  \label{fig:spy:rdb5000}
\end{subfigure}
\caption{Sparsity structures of the test matrices.}
\label{fig:spy}
\end{figure}

In this subsection we give the selection criterion and describe the matrices that we use for numerical experiments. {\tt Mps-rcm} makes incomplete LU based preconditioned iterative solvers very robust. Therefore,  to identify the most challenging problems, we use {\tt mps-rcm} to choose a highly indefinite and challenging set of $10$ problems from the SuiteSparse Matrix Collection~\cite{davis2011university} in which at least one instance of {\tt mps-rcm} fails due to failure of  incomplete factorization or stagnation of the Krylov subspace method. Properties of the test problems and their sparsity plots are given in Table~\ref{tab:properties} and Figure~\ref{fig:spy}, respectively. All chosen problems are (numerically) non-symmetric, and only a few of them are structurally symmetric. Bp\_200 and bp\_600 are from a sequence of simplex basis matrices in  Linear Programming. West0989 and west1505 arise in a chemical engineering plant model with seven and eleven stage column sections, respectively. Rajat19 is a circuit simulation problem. Rdb1250l, rdb3200l and rdb5000 arise in a reaction-diffusion Brusselator model.  Chebyshev2 is an integration matrix using the Chebyshev method for solving fourth-order semilinear initial boundary value problems and finally, Orani678 arises in the economic modeling of Australia.

\subsection{Effectiveness of the shifted skew-symmetrizer}

The structure of the approximate  skew-symmetrizer ($S$) can be anything. We experimented with a simple diagonal ($S_d$) and tridiagonal ($S_t)$ structures.
In Table~\ref{tab:lls_size}, the dimensions, and the number of nonzeros for the LLS problem in (\ref{lls_problem}) are given.  After that, we obtain a shifted skew-symmetrized matrix ($\widehat{A}$). To evaluate the effectiveness of the scaling and permutation followed by the approximate skew-symmetrizer, we use three metrics: the skew-symmetry of the off-diagonals, the distance of the main diagonal to identity, and the condition number. In Table \ref{tab:results}, we depict these for the original matrix, for matrices after MC64 scaling and permutation, and followed by applying  $S_d$ or $S_t$, which we call "Original", "MC64", "MC64+$S_d$" and "MC64+$S_t$", respectively. As expected, for most cases MC64+$S_t$ has successfully  improved the skew-symmetry of the off-diagonal part compared to the original matrix. One exception is chebyshev2, which has a condition number of order $10^{15}$ and MC64+$S_t$ has improved both the condition number and the main diagonal. For all test problems, MC64+$S_t$ has improved the main diagonal and for $8$  of $10$ cases, it has also improved the condition number compared to the original matrix.  In all cases, $S_t$ improves the skew-symmetry of the off-diagonal part and the diagonal compared to $S_d$  except chebyshev2. The condition number becomes worse for $6$ cases out of $10$ using $S_t$. However, this is not an issue since we further precondition the system in our proposed method.

\begin{figure}[ht]
\begin{subfigure}{.32\textwidth}
  \centering
  \includegraphics[width=1\linewidth]{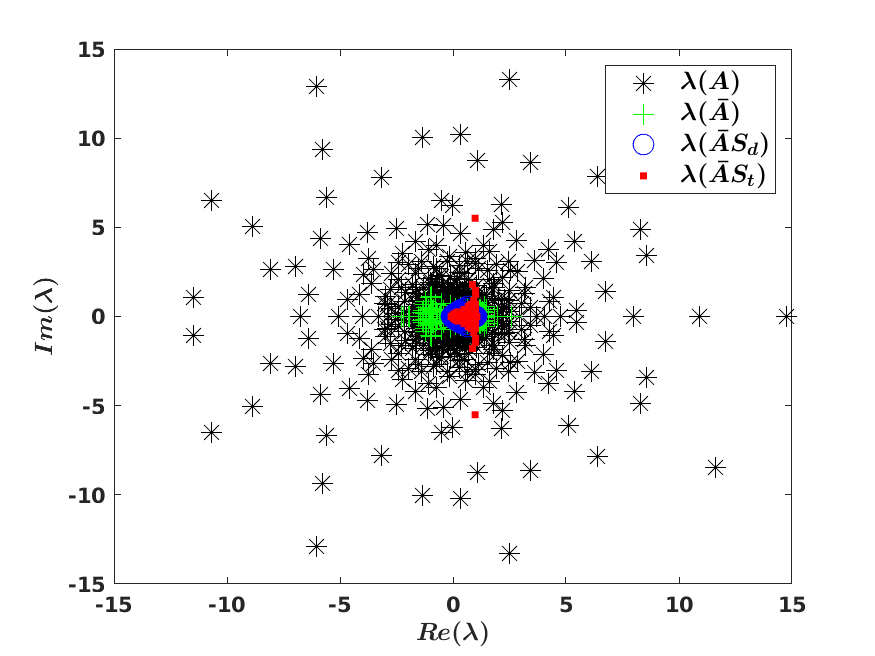}
  \caption{bp\_200}
  \label{fig:specs:bp_200}
\end{subfigure}%
\begin{subfigure}{.32\textwidth}
  \centering
  \includegraphics[width=1\linewidth]{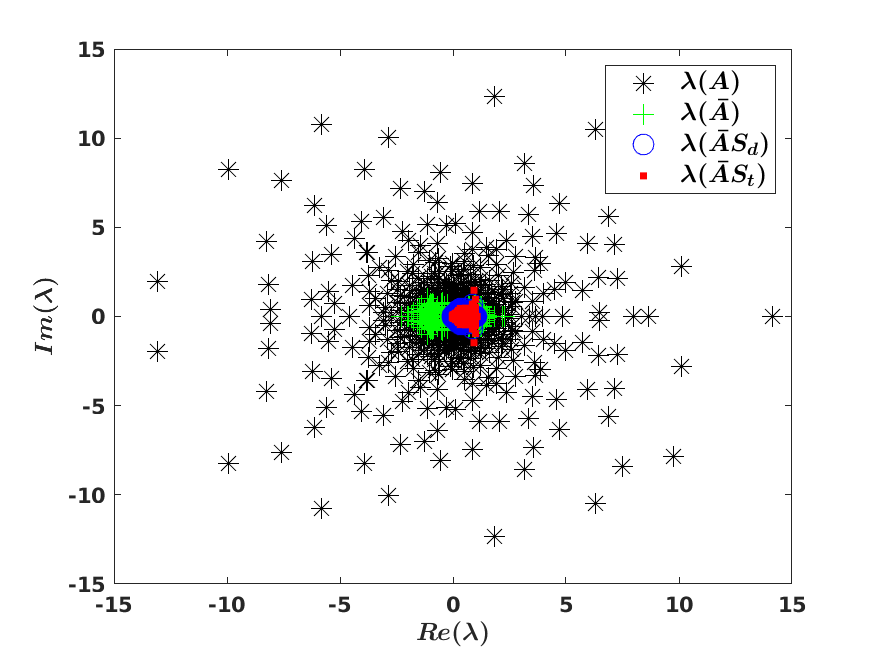}
  \caption{bp\_600}
  \label{fig:specs:bp_600}
\end{subfigure}
\begin{subfigure}{.32\textwidth}
  \centering
  \includegraphics[width=1\linewidth]{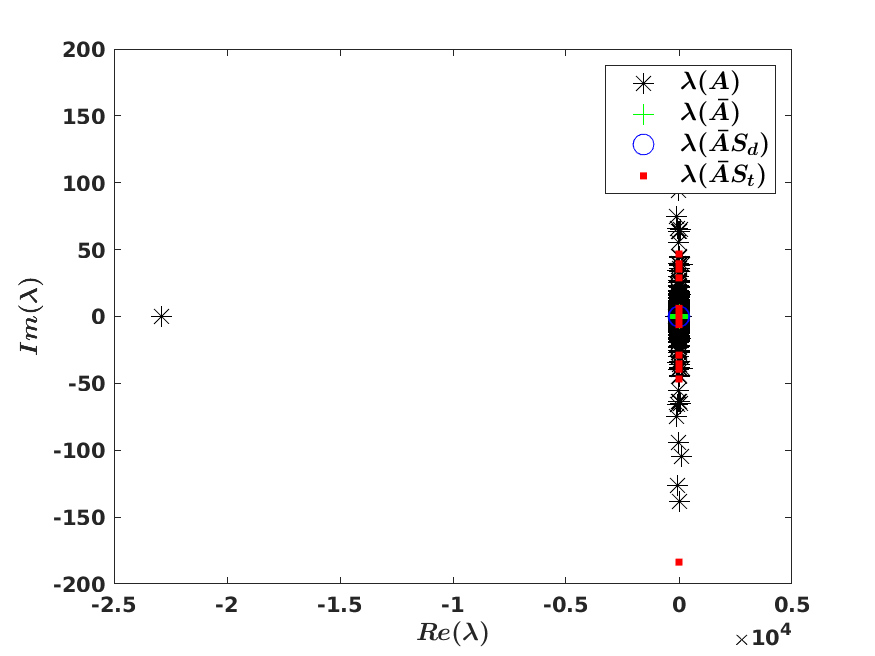}
  \caption{west0989}
  \label{fig:specs:west0989}
\end{subfigure}
\begin{subfigure}{.32\textwidth}
  \centering
  \includegraphics[width=1\linewidth]{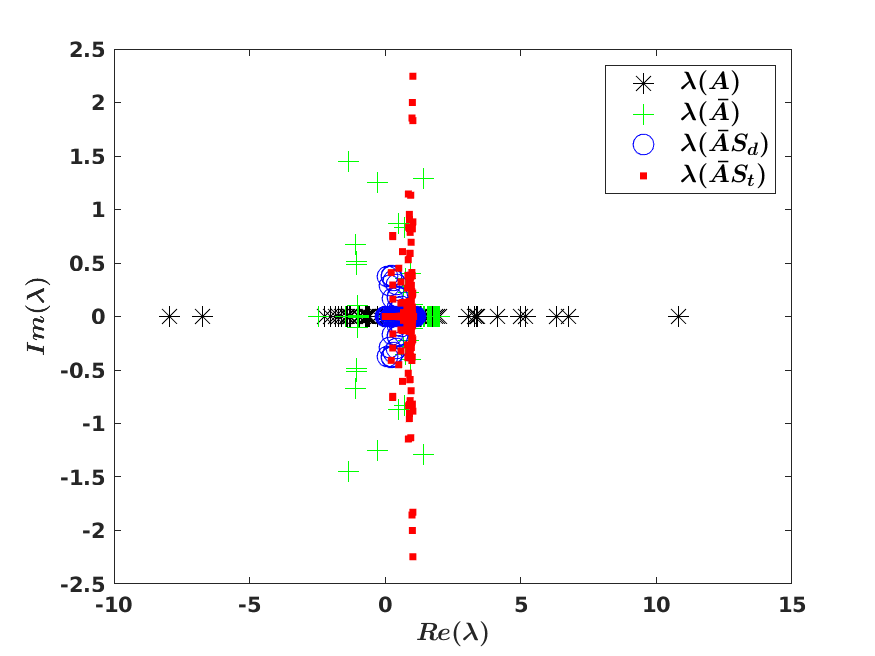}
  \caption{rajat19}
  \label{fig:specs:rajat19}
\end{subfigure}
\begin{subfigure}{.32\textwidth}
  \centering
  \includegraphics[width=1\linewidth]{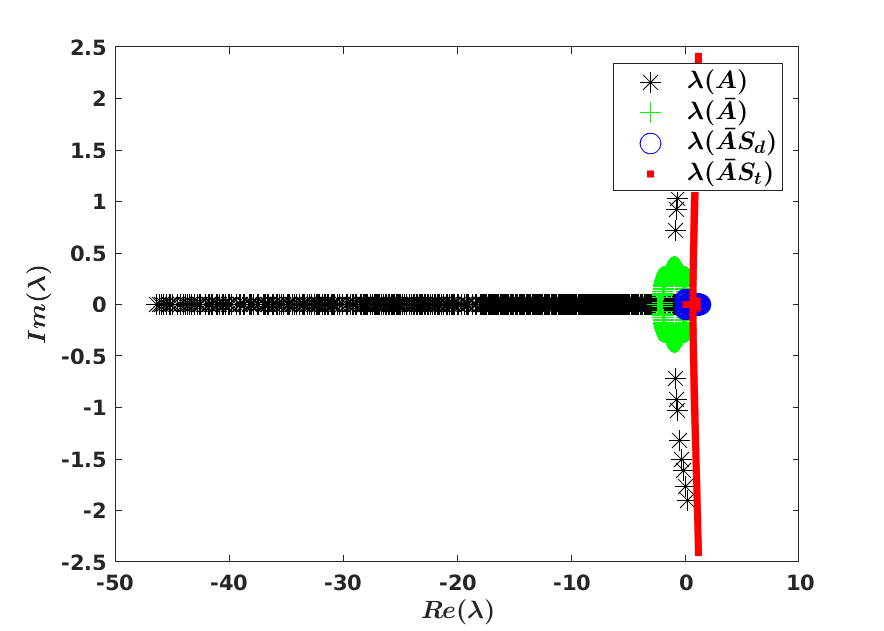}
  \caption{rdb1250l}
  \label{fig:specs:rdb1250l}
\end{subfigure}
\begin{subfigure}{.32\textwidth}
  \centering
  \includegraphics[width=1\linewidth]{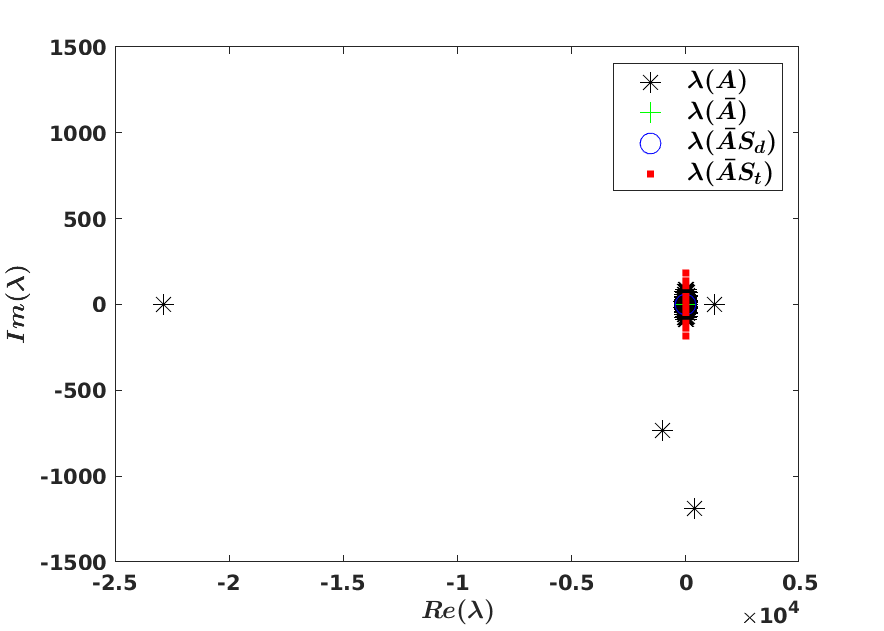}
  \caption{west1505}
  \label{fig:specs:west1505}
\end{subfigure}
\begin{subfigure}{.32\textwidth}
  \centering
  \includegraphics[width=1\linewidth]{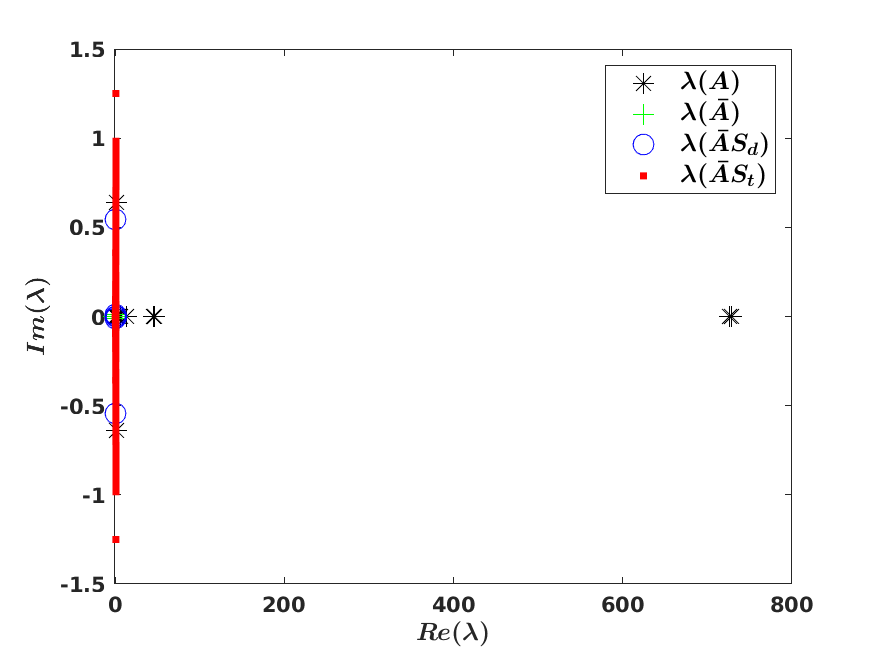}
  \caption{chebyshev2}
  \label{fig:specs:chebyshev2}
\end{subfigure}
\begin{subfigure}{.32\textwidth}
  \centering
  \includegraphics[width=1\linewidth]{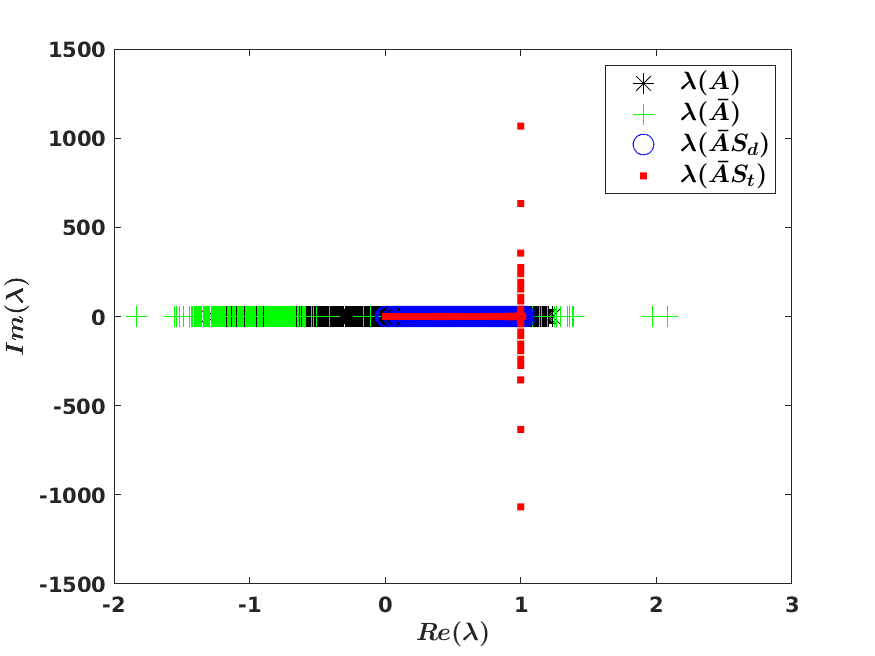}
  \caption{orani678}
  \label{fig:specs:orani678}
\end{subfigure}
\begin{subfigure}{.32\textwidth}
  \centering
  \includegraphics[width=1\linewidth]{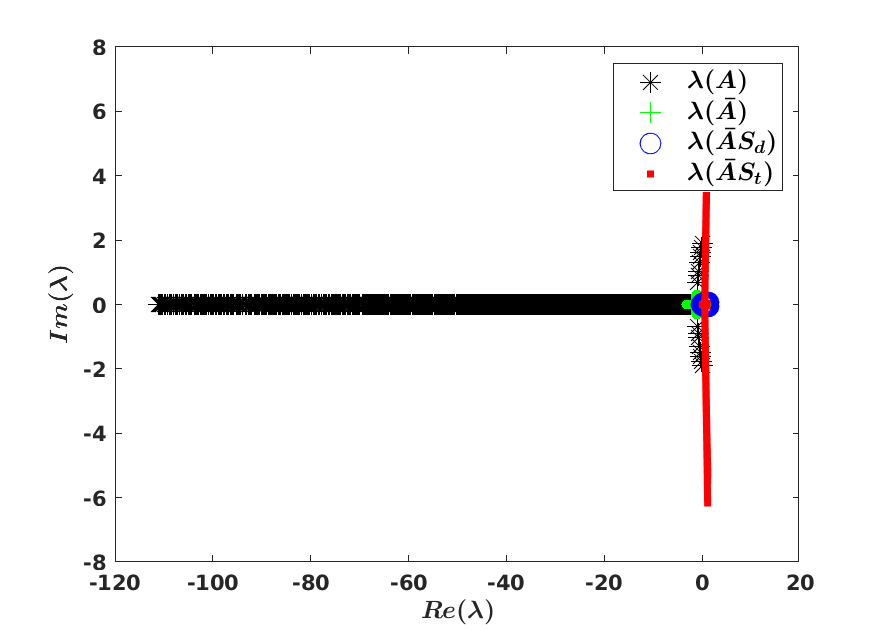}
  \caption{rdb3200l}
  \label{fig:specs:rdb3200l}
\end{subfigure}
\begin{subfigure}{\textwidth}
  \centering
  \includegraphics[width=0.32\linewidth]{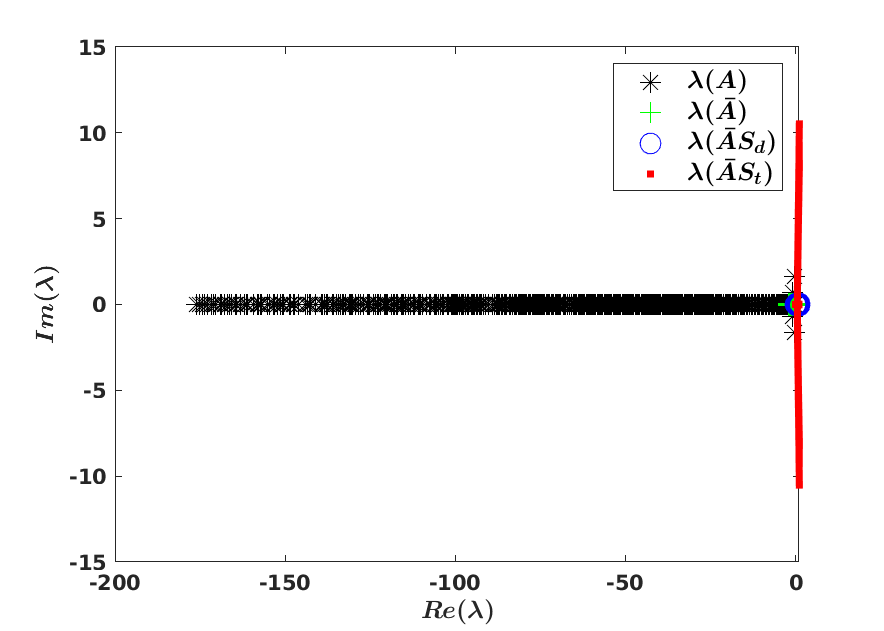}
  \caption{rdb5000}
  \label{fig:specs:rdb5000}
\end{subfigure}
\caption{Spectrum of the original matrix, after applying MC64 and shifted skew-symmetrizers.}
\label{fig:specs}
\end{figure}

Spectra of the original, reordered, and skew-symmetrized matrices are given in Figure~\ref{fig:specs}.  $S_t$ (shown in red) does a better job moving the real part of most eigenvalues positive side of the real axis and clustering them around one, compared to $S_d$.  Therefore, in the following numerical experiments we use $S_t$.

\begin{table}[ht]
\centering
\caption{Size (n), number of nonzeros (nnz), structural symmetry (Struct. S.), numerical symmetry (Num. S.) and the problem domains  of the test problems.}
\label{tab:properties}
\begin{tabular}{l | rrccl}
Matrix          & n   & nnz   & Struct. S. & Num. S. & Problem Domain               \\ \hline
bp\_200         & $822$  & $3,802$   & n    &   n       &Optimization\\
bp\_600         & $822$  & $4,172$   & n    & n     &Optimization\\
west0989        & $989$  & $3,518$   & n    & n       &Chemical Process Simulation\\
rajat19         & $1,157$ & $3,699$   & n    & n      &Circuit Simulation\\
rdb1250l        & $1,250$& $3,802$   & y    & n      &Computational Fluid dynamics\\
west1505        & $1,505$ & $5,414$  & n    & n      &Chemical Process Simulation\\
cehbyshev2      & $2,053$ & $18,447$ & n     & n     &Structural \\
orani678        & $2,529$ & $90,158$ & n     &n   &Economics \\
rdb3200l        & $3,200$  &$18,880$ & y    & n & Computational Fluid Dynamics \\
rdb5000         & $5,000$  &$29,600$& y     & n & Computational Fluid Dynamics
\end{tabular}
\end{table}
\begin{table}[ht]
\centering
\caption{Dimensions and number of nonzeros of the LLS problems associated with the test problems.}
\label{tab:lls_size}
\begin{tabular}{l|rrr|rrr}
                & \multicolumn{3}{c|}{$S_d$} & \multicolumn{3}{c}{$S_t$} \\
Matrix          & n       & m       & nnz      & n        & m        & nnz       \\\hline
bp\_200		& $3,722$ & $822$   & $3,802$  & $8,131$  & $2,463$    & $11,404$   \\
bp\_600		& $4,007$ & $822$   & $4,172$  & $9,604$  & $2,464$    & $12,514$   \\
west0989	& $3,491$ & $989$   & $3,518$  & $7,529$  & $2,965$    & $10,549$   \\
rajat19		& $3,425$ & $1,157$   & $3,699$  & $7,422$  & $3,469$    & $11,095$   \\
rdb1250l	& $4,275$ & $1,250$   & $7,300$  & $8,570$  & $3,748$    & $21,892$   \\
west1505	& $5,382$ & $1,505$   & $5,414$  & $11,598$  & $4,513$    & $16,237$   \\
chebyshev2& $14,344$ & $2,053$   & $18,447$  & $20,481$  & $6,157$    & $55,331$   \\
orani678	& $87,358$ & $2,529$   & $90,158$  & $112,852$  & $7,585$    & $270,465$   \\
rdb3200l & $11,040$ & $3,200$   & $18,880$  & $7,422$  & $3,469$  & $11,095$ \\
rdb5000 & $17,300$ & $5,000$   &  $29,600$ & $34,645$  & $14,998$    & $88,792$
\end{tabular}
\end{table}
\begin{table}[ht]
\centering
\caption{The effect of MC64 and the shifted skew symmetrizer.  In the table Skew-symmetry, Diagonal and Cond denote  $||(X-X^T)/2||_F/||X-\mathcal{D}(X)||_F$, $||\mathcal{D}(X)-I||_F$ and the condition number of $X$, respectively and rounded to one decimal place where $X$ is either the original matrix, after applying MC64, MC64 followed by $S_d$ or MC64 followed by $S_t$ in which $S_d$ and $S_t$ are diagonal and tridiagonal shifted skew-symetrizers, respectively.}
\label{tab:results}
\begin{tabular}{ll|llll}
                Matrix        &               & Original & MC64     & MC64+$S_d$  & MC64+$S_t$   \\\hline
\multirow{3}{*}{bp\_200}           & Skew-symmetry & 70.7\%   & 71\%   & 70.7\%   & 86.2\%    \\
                                 & Diagonal      & 28.7 & 35.8    & 12.2    &  11.6     \\
                                 & Cond          & $6.4\times 10^6$ & $4.0\times 10^2$ & $4.8\times 10^2$ & $1.0\times 10^4$  \\\hline
\multirow{3}{*}{bp\_600}        & Skew-symmetry & 70.7\%   &70.7\%   & 71.3\%   & 74.8\%    \\
                                 & Diagonal      & 28.7    & 36.6    & 13.2    & 12.7     \\
                                 & Cond          & $1.5\times 10^6$ & $3.2\times 10^2$ & $3.5\times 10^2$ & $7.8\times10^3$  \\\hline
\multirow{3}{*}{west0989}          & Skew-symmetry & 70.7\%   & 70.7\%   & 70.7\%   & 100\%    \\
                                 & Diagonal      & $2.3\times 10^{4}$    & 28.8     & 13.4     & 12.6      \\
                                 & Cond          & $9.9\times 10^{11}$ & $6.7\times 10^3$ & $8.3\times 10^3$ & $9.3\times 10^5$  \\\hline
\multirow{3}{*}{rajat19}         & Skew-symmetry & 28.4\%   & 67.1\%   & 68.1\%   & 83.4\%    \\
                                 & Diagonal      & 33.9   & 29.5    & 11.9    & 9.8      \\
                                 & Cond          & $1.1\times 10^{10}$ & $2.3\times 10^{10}$ & $1.1\times 10^{11}$ & $5.7\times 10^{10}$  \\\hline
\multirow{3}{*}{rdb1250l} & Skew-symmetry & 49\%    &56.2 \%   & 55.4\%   & 97.9\%    \\
                                 & Diagonal      & 690.8    & 70.7    & 17.2    & 9.8     \\
                                 & Cond          & $4.7\times 10^{2}$ & $4.9\times 10^2$ & $3.6\times 10^2$ & $3.1\times 10^2$  \\\hline
\multirow{3}{*}{west1505}        & Skew-symmetry & 70.7\%   & 70.7\%   & 70.7\%   & 100\%    \\
                                 & Diagonal      & $2.3\times 10^4$    & 35.9    & 16.6    & 15.7     \\
                                 & Cond          & $1.6\times 10^{12}$ & $8.8\times 10^3$ & $1.1\times 10^4$ & $1.2\times 10^6$\\\hline
\multirow{3}{*}{chebyshev2} & Skew-symmetry & 70.7\%    & 11.8\%   & 7.6\%   & 37.4\%    \\
                                 & Diagonal      & 898.5    & 3.5    & 21.9   & 21.9     \\
                                 & Cond          & $5.5\times 10^{15}$ & $8.6\times 10^9$ & $2.9\times 10^{10}$ & $1.5\times 10^{10}$  \\\hline
\multirow{3}{*}{orani678} & Skew-symmetry & 70.7\%    & 70.8\%   & 70.6\%   & 100\%    \\
                                 & Diagonal      & 53.5    & 97.5    & 15.2    & 15.1     \\
                                 & Cond          & $9.6\times 10^{3}$ & $7.5\times 10^3$ & $1.2\times 10^4$ & $6.4\times 10^6$  \\\hline
\multirow{3}{*}{rdb3200l} & Skew-symmetry & 21.9\%    & 27.9\%   & 27.1\%   & 99.7\%    \\
                                 & Diagonal      & $2.5\times 10^3$    & 113.1    & 20.6    & 12     \\
                                 & Cond          & $1.1\times 10^{3}$ & $9\times 10^2$ & $8.2\times 10^2$ & $7.3\times 10^2$  \\\hline
\multirow{3}{*}{rdb5000} & Skew-symmetry & 14.4\%    & 18.3\%   & 17.9\%   & 99.9\%    \\
                                 & Diagonal      & $4.8\times10^3$   & 141.4    & 24.6    & 14.1    \\
                                 & Cond          & $4.4\times 10^{3}$ & $3\times 10^3$ & $2.8\times 10^3$ & $3.6\times 10^3$  \\
\end{tabular}
\end{table}
\begin{table}[ht]
\centering
\caption{Rank of $\mathcal{\tilde{M}}_r$ and its percentage with respect to the matrix dimension ($\frac{r}{n}\times 100$) rounded to one decimal place in parenthesis.}
\label{tab:rank}
\begin{tabular}{l|lll}
Matrix          & $\sim ildl(0)$ & $ildl(10^{-1})$       & $ildl(10^{-2})$         \\\hline
bp\_200		& $33(4\%)$ & $32(3.8\%)$   & $32(3.8\%)$    \\
bp\_600		& $48(5.8\%)$ & $46(5.6\%)$   & $45(5.5\%)$     \\
west0989	& $35(3.5\%)$ & $32(3.2\%)$   & $35(3.5\%)$    \\
rajat19		& $38(3.3\%)$ & $38(3.3\%)$   & $38(3.3\%)$     \\
rdb1250l	& $36(2.9\%)$ & $1(0.1\%)$   & $26(2.1\%)$     \\
west1505	& $54(3.6\%)$ & $55(3.7\%)$   & $56(3.7\%)$    \\
chebyshev2& $5(2\%)$ & $7(0.3\%)$   & $7(0.3\%)$    \\
orani678	& $19(0.8\%)$ & $26(1\%)$   & $28(1.1\%)$     \\
rdb3200l & $0(0\%)$ & $0(0\%)$   & $1(0\%)$   \\
rdb5000 & $0(0\%)$ & $0(0\%)$   &  $0(0\%)$
\end{tabular}
\end{table}

\subsection{Effectiveness of the deflation}\label{sec:deflation}
\begin{figure}[ht]
\begin{subfigure}{.32\textwidth}
  \centering
  \includegraphics[width=1\linewidth]{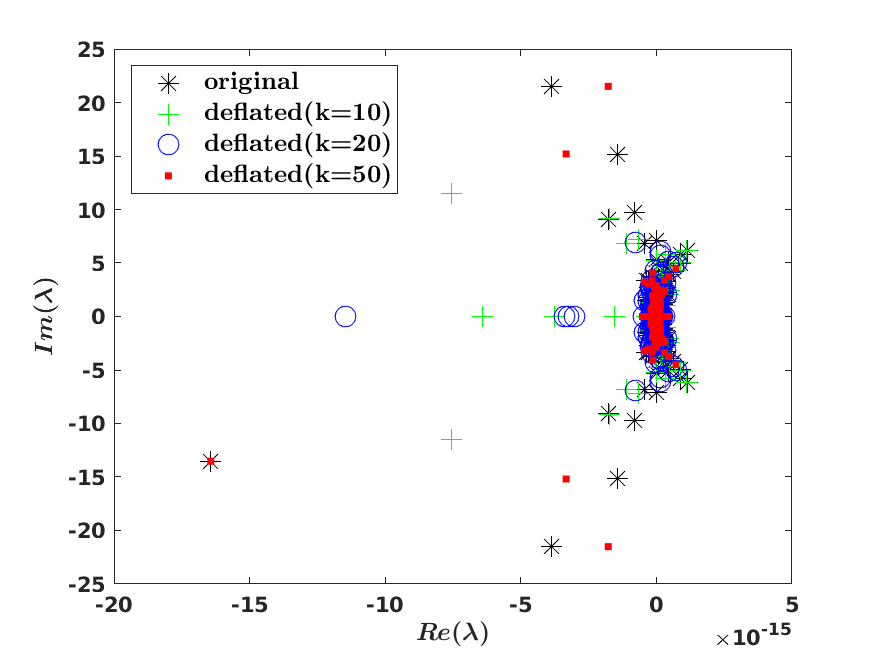}
  \caption{bp\_200}
  \label{fig:def:bp_200}
\end{subfigure}%
\begin{subfigure}{.32\textwidth}
  \centering
  \includegraphics[width=1\linewidth]{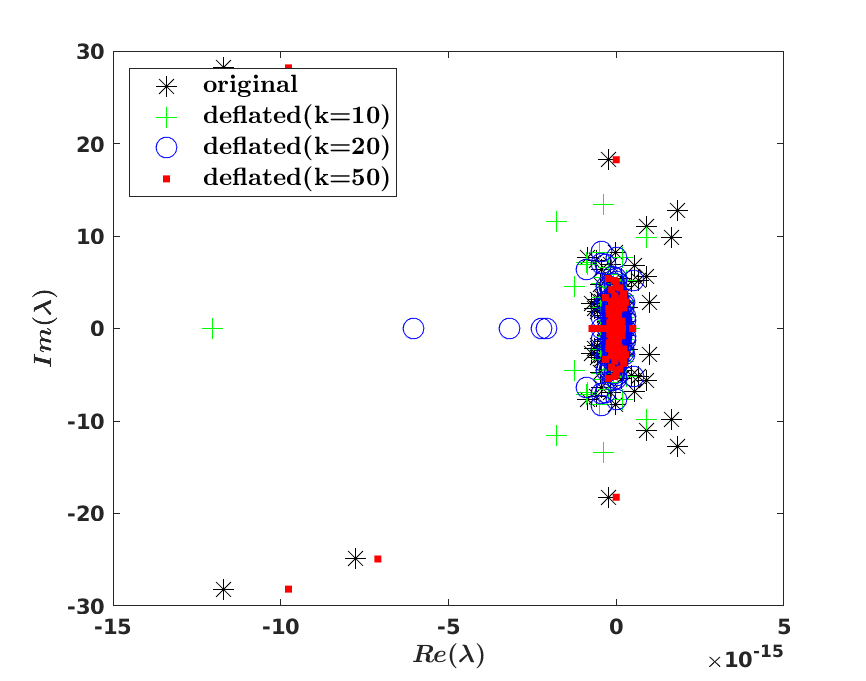}
  \caption{bp\_600}
  \label{fig:def:bp_600}
\end{subfigure}
\begin{subfigure}{.32\textwidth}
  \centering
  \includegraphics[width=1\linewidth]{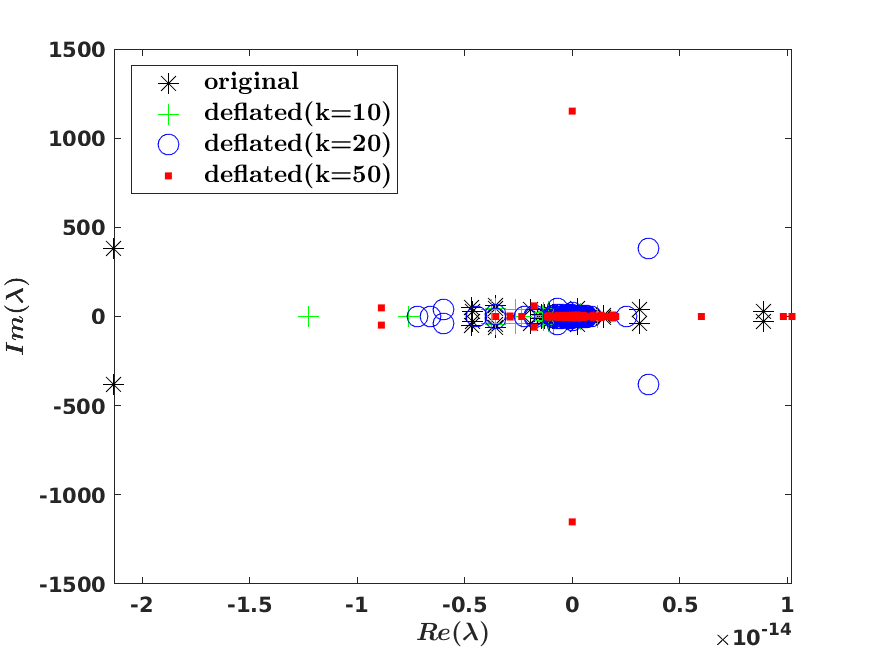}
  \caption{west0989}
  \label{fig:def:west0989}
\end{subfigure}
\begin{subfigure}{.32\textwidth}
  \centering
  \includegraphics[width=1\linewidth]{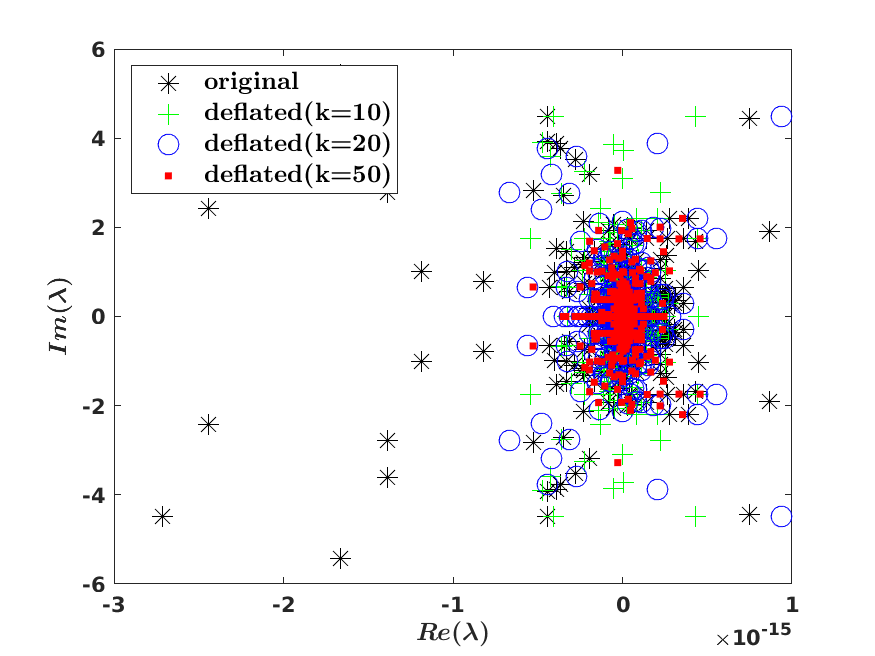}
  \caption{rajat19}
  \label{fig:def:rajat19}
\end{subfigure}
\begin{subfigure}{.32\textwidth}
  \centering
  \includegraphics[width=1\linewidth]{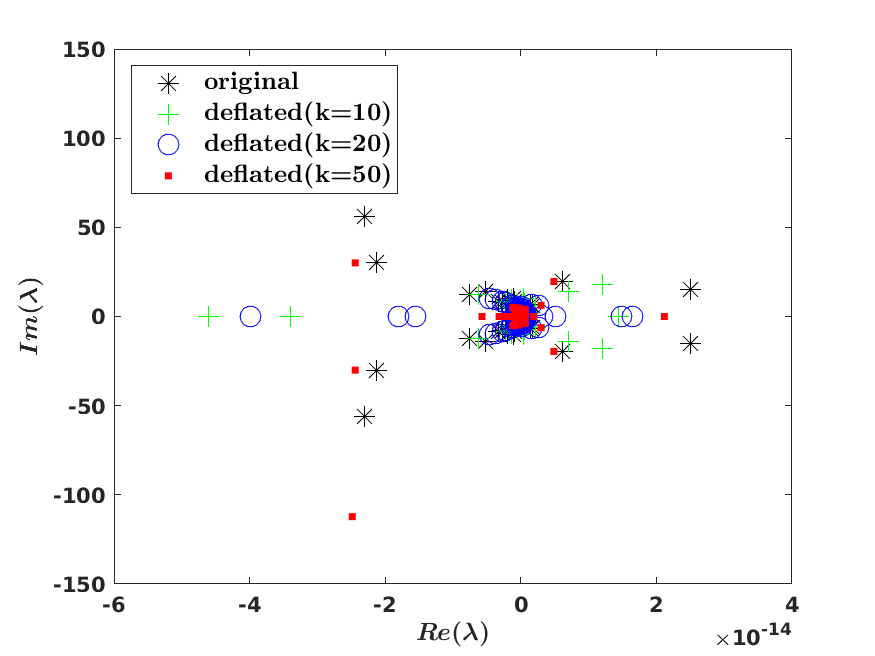}
  \caption{rdb1250l}
  \label{fig:def:rdb1250l}
\end{subfigure}
\begin{subfigure}{.32\textwidth}
  \centering
  \includegraphics[width=1\linewidth]{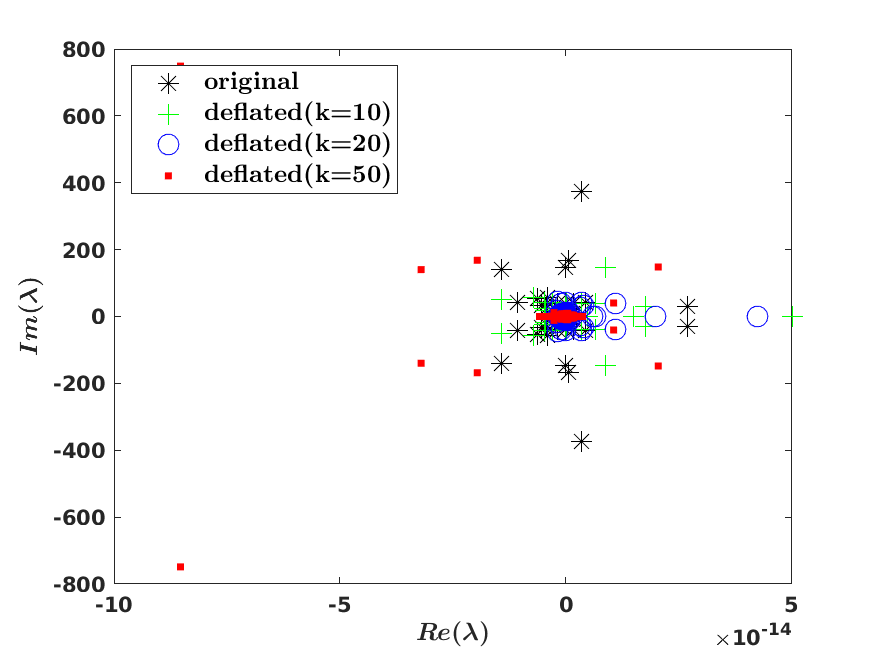}
  \caption{west1505}
  \label{fig:def:west1505}
\end{subfigure}
\begin{subfigure}{.32\textwidth}
  \centering
  \includegraphics[width=1\linewidth]{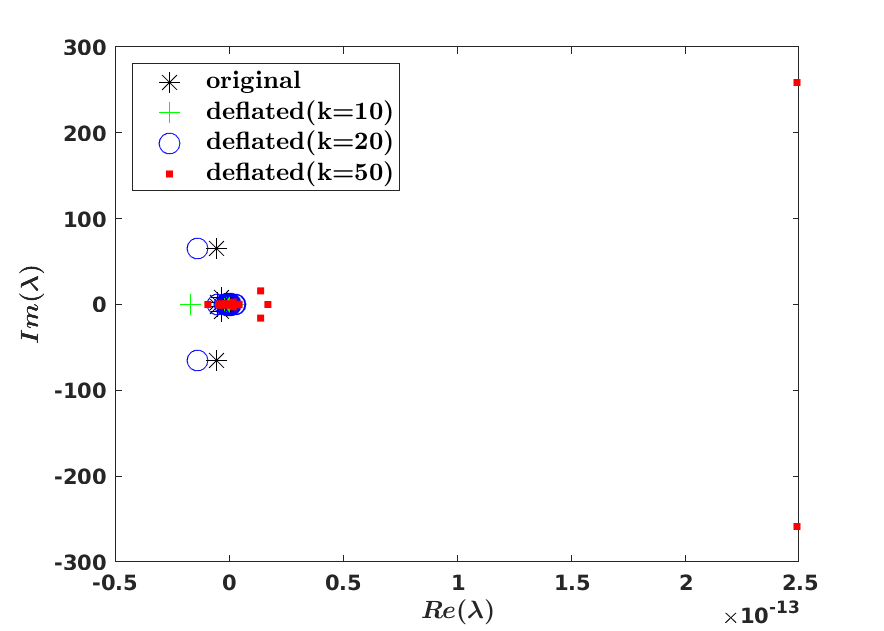}
  \caption{chebyshev2}
  \label{fig:def:chebyshev2}
\end{subfigure}
\begin{subfigure}{.32\textwidth}
  \centering
  \includegraphics[width=1\linewidth]{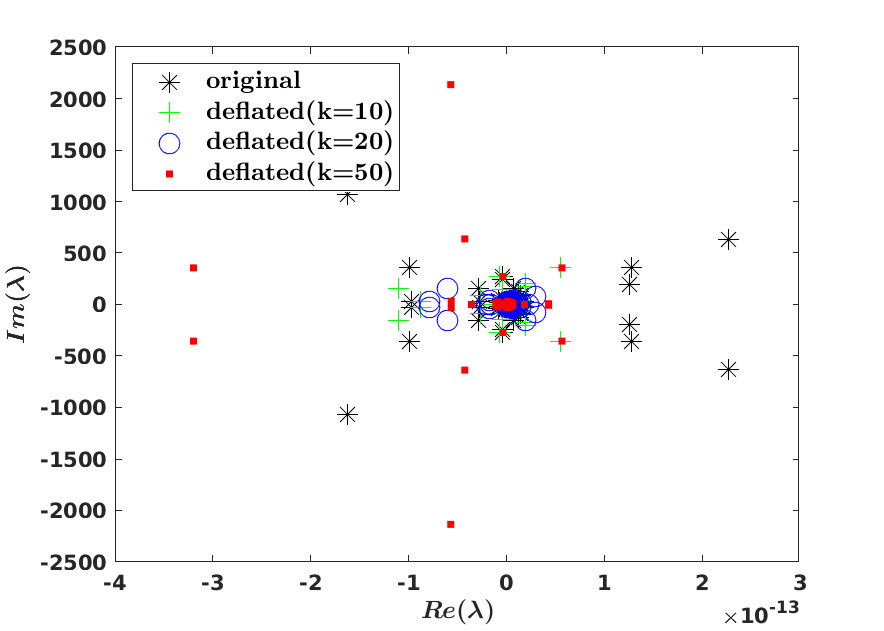}
  \caption{orani678}
  \label{fig:def:orani678}
\end{subfigure}
\begin{subfigure}{.32\textwidth}
  \centering
  \includegraphics[width=1\linewidth]{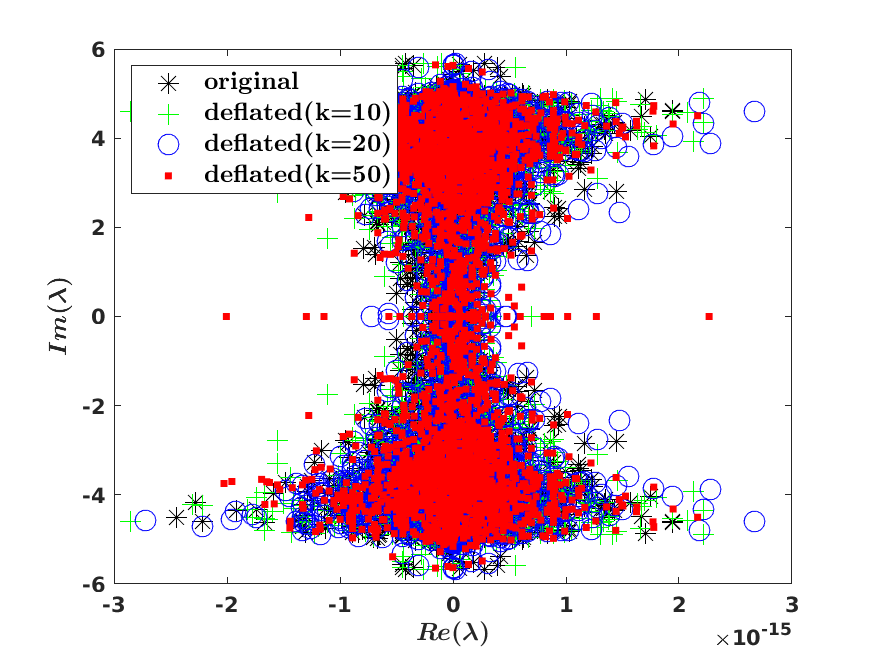}
  \caption{rdb3200l}
  \label{fig:def:rdb3200l}
\end{subfigure}
\begin{subfigure}{\textwidth}
  \centering
  \includegraphics[width=0.32\linewidth]{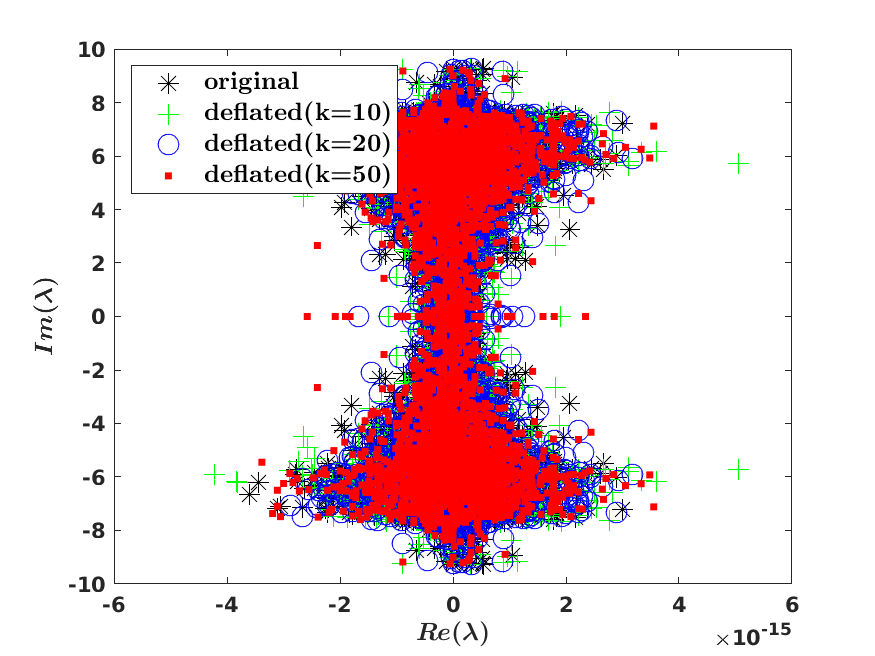}
  \caption{rdb5000}
  \label{fig:def:rdb5000}
\end{subfigure}
\caption{Spectrum of the original skew-symmetric matrix and after deflation.}
\label{fig:def}
\end{figure}

In Figure~\ref{fig:def}, we present the spectrum of the original $I+\tilde{\mathcal{J}}$ (by using incomplete $LDL^{T}$ factorization of $\widehat{H}$ with zero fill in) and of the deflated shifted skew-symmetric $I+\bar{\mathcal{J}}$ after $10$, $20$ and $50$ iterations of the  skew-Lanczos process for all test matrices. We note that the scale of the real axis is negligible for all cases and the spectrum is  purely imaginary. Even though case by case fine-tuning is possible by looking at the spectra, we observe that deflating with $20$ vectors gives a meaningful balance between the number of iterations and the orthogonality between the Lanczos vectors for the test problems, since for $50$ vectors the spectrum is worse  and for $10$ vectors it does not improve it as much as for $20$ vectors. Therefore, in the following experiments, we use $20$ skew-Lanczos vectors.
\subsection{Iterative solution of general sparse linear systems}
In the following numerical experiments we use the proposed method as described earlier. For the proposed method the number of Lanczos vectors is set to $20$. Compared to the non-deflated version of the proposed method, the number of {\tt mrs} iterations are $30.8\%$ better on average when the inner system is deflated with $20$ Lanczos vectors for the test problems. We have also experimented with the skew-Lanczos method with full reorthogonalization for $20$ Lanczos vectors and the improvement in the number of iterations is negligible for the test problems. Since it is used just as a preconditioner, we believe that in general it may not be necessary to fully or partially reorthgonalize in the skew-Lanczos process. Therefore, we only report the proposed method with deflation using $20$ Lanczos vectors without any reorthogonalization.

In Table \ref{tab:rank}, the ranks  of $\tilde{\mathcal{M}}_r$ for various incomplete factorizations, as well as the percentages of the rank with respect to the matrix dimension are given. For the set of test problems the largest one is $56$, which is only $3.7\%$ of the matrix dimension. They are roughly invariant for incomplete factorization for all problems except $rdbl1250l$.  For the cases when the rank is zero, we use the preconditioner in (\ref{preconditioner2}), otherwise we use the preconditioner in (\ref{preconditioner}). We have also experimented with the preconditioner in (\ref{preconditioner2}) for problems where the rank is not zero but relatively small. However, the number of iterations has increased significantly even for the case where rank is equal to one ($rdbl3200l$ with $ildl(10^{-2})$). We have tried to solve the LLS problem with the regularization parameters (\ref{eq:regularized}) $\gamma=0.1, 1$ and $10$. Since $\gamma=0.1$ and $10$ are worse in terms of the number of iterations, we only report the results with $\gamma=1$.

In Table \ref{tab:iterations}, the number of {\tt tfqmr} iterations for the proposed method and {\tt mps-rcm} are given.  For the proposed method, we also give the average number of inner iterations in parenthesis. As seen in the table, for all test problems and regardless of the choice of the dropping tolerance, the proposed method succeeds. On the other hand, {\tt mps-rcm} fails for $80\%$, $90\%$, and $20\%$ of problems using $ilu(0)$, $ilutp(10^{-1})$, and $ilutp(10^{-2})$, respectively. For $ilu(0)$ the majority of the failures are due to {\tt tfqmr} stagnating. While for $ilutp(10^{-1})$ it happens because of a zero pivot encountered during factorization. As expected,  when the dropping tolerance  decreases to $10^{-2}$ the failures also decrease since the incomplete factorization is more like a direct solver. Hence, {\tt mps-rcm} becomes more robust.  When it does not fail, the required number of {\tt tfqmr} iterations for {\tt mps-rcm} is quite low except for $rdb5000$.
The proposed method, on the other hand, is robust regardless of the quality of the incomplete factorization for the test problems. There are no failures during the incomplete  factorization and no failures of the iterative scheme. The required number of {\tt tfqmr} iterations improve as a more accurate incomplete factorization is used. The number of iterations if $ildl(10^{-2})$ is used are comparable to those obtained by {\tt mps-rcm} with $ilutp(10^{-2})$. Except,  for two cases ($rdb3200l$ and $rdb5000$) for which the proposed method is significantly better and for two other cases ($orani678$ and $rdb1250l$) for which {\tt mps-rcm} is significantly better. For the proposed method the number of average inner {\tt mrs} iterations is not dependent on the choice  of the incomplete factorization for all test problems except for $rdb1250l$; in this case, it is almost halved when $ildl(10^{-1})$ is used compared to $\sim ildl(0)$ and $ildl(10^{-2})$. We believe this is due to the fact that for the same matrix using $ildl(10^{-1})$, the rank of $\tilde{\mathcal{M}}_r$ is $1$ which is much smaller than those of $\sim ildl(0)$ and $ildl(10^{-1})$, $36$ and $26$, respectively.
\begin{table}[ht]
\centering
\caption{Number of {\tt tfqmr} iterations. The average number of inner {\tt mrs} iterations for the proposed method rounded to one decimal place is given in parenthesis. $\ast$: {\tt tfqmr} stagnated, $\dagger$:
{\tt tfqmr} reached the maximum number of iterations ($2,000$) without reaching the required relative residual, $\ddagger$: zero pivot is encountered during factorization. }
\label{tab:iterations}
\begin{tabular}{l|lll|lll}
                & \multicolumn{3}{c|}{{\bf Proposed method}} & \multicolumn{3}{c}{\bf mps-rcm} \\
Matrix          & $\sim ildl(0)$ & $ildl(10^{-1})$       & $ildl(10^{-2})$    &  $ilu(0)$       & $ilutp(10^{-1})$        & $ilutp(10^{-2})$       \\\hline
bp\_200		& $26(53.6)$ & $88(54.4)$   & $3(55.8)$  & $\dagger$  & $\ddagger$    & $2$   \\
bp\_600		& $15(65.7)$ & $33(66.6)$   & $4(67)$  & $\dagger$  & $\ddagger$    & $2$   \\
west0989	& $74(134.6)$ & $45(106.9)$   & $3(105.5)$  & $\ast$  & $20$    & $1$   \\
rajat19		& $10(35.6)$ & $46(35.8)$   & $17(35.2)$  & $\ddagger$  & $\ddagger$    & $\ddagger$   \\
rdb1250l	& $90(66.7)$ & $10(37.1)$   & $49(71.8)$  & $\ast$  & $\ddagger$    & $10$   \\
west1505	& $57(162.9)$ & $75(166.9)$   & $1(163)$  & $\dagger$  & $\ddagger$    & $1$   \\
chebyshev2& $10(15.9)$ & $15(12)$   & $14(15.6)$  & $1$  & $\ddagger$    & $\ddagger$   \\
orani678	& $15(229.3)$ & $13(235.9)$   & $12(249)$  & $15$  & $\ddagger$    & $4$   \\
rdb3200l & $9(63.6)$ & $15(61)$   & $3(69.1)$  & $\ast$  & $\ddagger$  & $20$ \\
rdb5000 & $8(100.6)$ & $16(101.5)$   &  $6(100.9)$ & $\ast$  & $\dagger$    & $129$
\end{tabular}
\end{table}
\section{Conclusions}
\label{sec:conclusions}
A robust two-level iterative scheme is presented for solving general sparse linear system of equations. The robustness of the scheme is shown on challenging matrices that arise in various problems which are obtained from the SuiteSparse Matrix Collection. While it requires some additional preprocessing steps, the results presented in this paper indicate that the proposed scheme significantly improves the robustness of iterative methods for general sparse linear systems compared to existing methods. This result is irrespective of the quality of the incomplete factorization, even for a challenging set of test problems. The proposed scheme requires some additional memory, but this is shown to be kept within a small percentage of the problem size. We believe with the introduction of the proposed scheme, iterative solvers will become  much more viable alternatives for solving  problems  in chemical engineering, optimizations, economics, etc. Its efficient parallel implementation requires addressing some algorithmic challenges which we leave as future work.


\bibliographystyle{siamplain}


\begin{thebibliography}{10}

\bibitem{alappat2020recursive}
{\sc C.~Alappat, A.~Basermann, A.~R. Bishop, H.~Fehske, G.~Hager, O.~Schenk,
  J.~Thies, and G.~Wellein}, {\em A recursive algebraic coloring technique for
  hardware-efficient symmetric sparse matrix-vector multiplication}, ACM
  Transactions on Parallel Computing (TOPC), 7 (2020), pp.~1--37.

\bibitem{BeaMXZ18}
{\sc C.~Beattie, V.~Mehrmann, H.~Xu, and H.~Zwart}, {\em Port-{H}amiltonian
  descriptor systems}, Math. Control Signals Systems, 30:17 (2018), p.~Appeared
  electronically.
\newblock https://doi.org/10.1007/s00498-018-0223-3.

\bibitem{BenHT00}
{\sc M.~Benzi, J.~C. Haws, and M.~Tuma}, {\em Preconditioning highly indefinite
  and nonsymmetric matrices}, SIAM Journal on Scientific Computing, 22 (2000),
  pp.~1333--1353.

\bibitem{BunKP76}
{\sc J.~R. Bunch, L.~Kaufman, and B.~N. Parlett}, {\em Decomposition of a
  symmetric matrix}, Numerische Mathematik, 27 (1976), pp.~95--109.

\bibitem{concus1976generalized}
{\sc P.~Concus and G.~H. Golub}, {\em A generalized conjugate gradient method
  for nonsymmetric systems of linear equations}, in Computing Methods in
  Applied sciences and Engineering, Springer, 1976, pp.~56--65.

\bibitem{ccuugu2020parallel}
{\sc {\.I}.~{\c{C}}u{\u{g}}u and M.~Manguo{\u{g}}lu}, {\em A parallel
  multithreaded sparse triangular linear system solver}, Computers \&
  Mathematics with Applications, 80 (2020), pp.~371--385.

\bibitem{davis2011university}
{\sc T.~A. Davis and Y.~Hu}, {\em The university of florida sparse matrix
  collection}, ACM Transactions on Mathematical Software (TOMS), 38 (2011),
  pp.~1--25.

\bibitem{DopU16}
{\sc F.~Dopico and F.~Uhlig}, {\em Computing matrix symmetrizers, part 2: New
  methods using eigendata and linear means; a comparison}, Linear Algebra and
  its Applications, 504 (2016), pp.~590--622.

\bibitem{DufK01}
{\sc I.~S. Duff and J.~Koster}, {\em On algorithms for permuting large entries
  to the diagonal of a sparse matrix}, SIAM Journal on Matrix Analysis and
  Applications, 22 (2001), pp.~973--996.

\bibitem{Egg19}
{\sc H.~Egger}, {\em Structure preserving approximation of dissipative
  evolution problems}, Numerische Mathematik, 143 (2019), pp.~85--106.

\bibitem{Fre93}
{\sc R.~W. Freund}, {\em A transpose-free quasi-minimal residual algorithm for
  non-hermitian linear systems}, SIAM Journal on Scientific Computing, 14
  (1993), pp.~470--482.

\bibitem{FreN91}
{\sc R.~W. Freund and N.~M. Nachtigal}, {\em {QMR}: a quasi-minimal residual
  method for non-hermitian linear systems}, Numerische Mathematik, 60 (1991),
  pp.~315--339.

\bibitem{GolV96}
{\sc G.~H. Golub and C.~F. {Van Loan}}, {\em Matrix Computations}, The Johns
  Hopkins University Press, Baltimore, MD, 3rd~ed., 1996.

\bibitem{GraMQSW16}
{\sc N.~Gr{\"a}bner, V.~Mehrmann, S.~Quraishi, C.~Schr\"oder, and U.~{von
  W}agner}, {\em Numerical methods for parametric model reduction in the
  simulation of disc brake squeal}, Z. Angew. Math. Mech., 96 (2016),
  pp.~1388--1405.

\bibitem{greif2017sym}
{\sc C.~Greif, S.~He, and P.~Liu}, {\em Sym-ildl: Incomplete ldlt factorization
  of symmetric indefinite and skew-symmetric matrices}, ACM Transactions on
  Mathematical Software (TOMS), 44 (2017), pp.~1--21.

\bibitem{GreV09}
{\sc C.~Greif and J.~M. Varah}, {\em Iterative solution of skew-symmetric
  linear systems}, SIAM Journal on Matrix Analysis and Applications, 31 (2009),
  pp.~584--601.

\bibitem{hsl2007collection}
{\sc HSL}, {\em A collection of fortran codes for large scale scientific
  computation.}, See \url{http://www.hsl.rl.ac.uk},  (2007).

\bibitem{IdeV07}
{\sc R.~Idema and C.~Vuik}, {\em A minimal residual method for shifted
  skew-symmetric systems}, Tech. Report REPORT 07-09, Delft University of
  Technology, 2007.

\bibitem{JacZ12}
{\sc B.~Jacob and H.~Zwart}, {\em Linear port-{H}amiltonian systems on
  infinite-dimensional spaces}, Operator Theory: Advances and Applications,
  223, Birkh{\"a}user/Springer Basel AG, Basel CH, 2012.

\bibitem{ji2017breakdown}
{\sc H.~Ji and Y.~Li}, {\em A breakdown-free block conjugate gradient method},
  BIT Numerical Mathematics, 57 (2017), pp.~379--403.

\bibitem{Jia07}
{\sc E.~Jiang}, {\em Algorithm for solving shifted skew-symmetric linear
  system}, Frontiers of Mathematics in China, 2 (2007), pp.~227--242.

\bibitem{manguoglu2019robust}
{\sc M.~Manguoglu and V.~Mehrmann}, {\em A robust iterative scheme for
  symmetric indefinite systems}, SIAM Journal on Scientific Computing, 41
  (2019), pp.~A1733--A1752.

\bibitem{PorU16}
{\sc R.~Portase and B.~U{\c c}ar}, {\em {On matrix symmetrization and sparse
  direct solvers}}, Research Report RR-8977, {Inria - Research Centre Grenoble
  -- Rh{\^o}ne-Alpes}, Nov. 2016, \url{https://hal.inria.fr/hal-01398951}.

\bibitem{SaaS86}
{\sc Y.~Saad and M.~H. Schultz}, {\em {GMRES}: A generalized minimal residual
  algorithm for solving nonsymmetric linear systems}, SIAM Journal on
  Scientific and Statistical Computing, 7 (1986), pp.~856--869.

\bibitem{Sch13}
{\sc A.~J.~v. Schaft}, {\em Port-{H}amiltonian differential-algebraic systems},
  in Surveys in Differential-Algebraic Equations I, Springer, 2013,
  pp.~173--226.

\bibitem{SzyW93}
{\sc D.~B. Szyld and O.~B. Widlund}, {\em Variational analysis of some
  conjugate gradient methods}, East-West Journal of Numerical Mathematics, 1
  (1993), pp.~51--74.

\bibitem{Uce08}
{\sc B.~U{\c c}ar}, {\em Heurisstics for a matrix symmetrization problem}, in
  Parallel Processing and Applied Mathematics, R.~Wyrzykowski, J.~Dongarra,
  K.~Karczewski, and J.~Wasniewski, eds., Berlin, Heidelberg, 2008, Springer
  Berlin Heidelberg, p.~718–727.

\bibitem{SchJ14}
{\sc A.~J. {van der}~Schaft and D.~Jeltsema}, {\em Port-{H}amiltonian systems
  theory: An introductory overview}, Foundations and Trends in Systems and
  Control, 1 (2014), pp.~173--378.

\bibitem{Van92}
{\sc H.~A. Van~der Vorst}, {\em {Bi-CGSTAB}: A fast and smoothly converging
  variant of bi-cg for the solution of nonsymmetric linear systems}, SIAM
  Journal on Scientific and Statistical Computing, 13 (1992), pp.~631--644.

\bibitem{vecharynski2013absolute}
{\sc E.~Vecharynski and A.~V. Knyazev}, {\em Absolute value preconditioning for
  symmetric indefinite linear systems}, SIAM Journal on Scientific Computing,
  35 (2013), pp.~A696--A718.

\bibitem{Wid78}
{\sc O.~Widlund}, {\em A {L}anczos method for a class of nonsymmetric systems
  of linear equations}, SIAM Journal on Numerical Analysis, 15 (1978),
  pp.~801--812.

\bibitem{zhang2000preconditioned}
{\sc J.~Zhang}, {\em Preconditioned {K}rylov subspace methods for solving
  nonsymmetric matrices from cfd applications}, Computer methods in Applied
  Mechanics and Engineering, 189 (2000), pp.~825--840.

\end{thebibliography}
\appendix
\section{Shifted skew-symmetrizer matlab code} \label{app:A}

\begin{PseudoCode}
function S = skew_symmetrize_sparse_right(A,n,gamma,S)
% Find S such that off-diaognals of AS are skew-symmetric
% and the diagonal of AS is as close to the identity as possible.
%
% Input
%   A: square (nxn), sparse input matrix
%   n: size of A
%   gamma: a positive regularization parameter
%   S: a sparsity pattern of the S matrix
% Output
%   S: the given sparsity structure such that AS is as
%      shifted skew-symetrix as possible
%

% compute the sparsity structure of AS
AA = spones(A)*spones(S);

% sparsity structure of the strictly upper triangular part
% each nonzero defines one equality in the LLS problem
U = triu(AA+AA',1);

% number of equations of the LLS  for the symmetry requirement only
neq = nnz(U);

% and their indices
[I,J] = find(U);

% extract the indices of unknowns in column major order
II = find(S);

% number of unknowns
nu = length(II);

% place "the compressed" indices in the corresponding nonzero entries
S(II) = [1:nu];

% Allocate the LLS coeffcient matrix
B = sparse(neq+n,nu);

% for each nonzero: (AS)ij
for i=1:neq
    % Our equality is:
    %       (AS)ij + (AS)ji = 0
    % first we need to determine those indices of nonzeros that contribute to the equality
    ASij = spones(A(I(i),:))'.*spones(S(:,J(i)));
    ASji = spones(A(J(i),:))'.*spones(S(:,I(i)));

    INDij=find(ASij);
    B(i,S(INDij,J(i))) = A(I(i),INDij);

    INDji=find(ASji);
    B(i,S(INDji,I(i))) = A(J(i),INDji);

end

% additional constraints to obtain ones on the main diagonal
for i=1:n
   ASii =   spones(A(i,:))'.*spones(S(:,i));
   INDii=find(ASii);
   B(i+neq,S(INDii,i)) = sqrt(gamma)*A(i,INDii);
end

% set the right hand side vector
v = [zeros(neq,1);sqrt(gamma)*ones(n,1)];

% Solve the sparse LLS problem
x = B\v;

% map solution back to S matrix
S(II) = x;

return
\end{PseudoCode}

\section{Shifted skew-symmetric iterative solver for multiple right hand side vectors based on simultaneous {\em mrs} iterations} \label{app:B}

\begin{PseudoCode}
function [x, its, res, relres_hist] = mrs(alpha,S,b,maxit,tol)
%solves linear systems (alpha*I + S )x = b
% where b and x are of size nxnrhs
% Input
%   alpha   : a scalar
%   S       : a sparse skew-symmetrix matrix
%   b       : the right hand side vector
%   maxit   : maximum number of iteratations
%   tol     : stopping tolerance
% Output
%   x           : solution vector
%   its         : number of iterations
%   res         : final residuals
%   relres_hist : relative residual history

[n,nrhs] = size(b)

x = zeros(n,nrhs);
r = b;

for j =1:nrhs
    r0(j) = norm(r(:,j),2);
    relres_hist(j) = [r0(j)/r0(j)];
    s(j) = r0(j);
    q(:,j) = r(:,j)/s(j);
    beta(j) = 0;
    theta1(j) = alpha;
    c_old(j) = 1;
    s_old(j) = 0;
    delta(j) = 0;
    delta_old(j) = 0;
end
q_old = zeros(n,nrhs);
p_old = zeros(n,nrhs);
p_old2 = zeros(n,nrhs);


for i=1:maxit
    q_new = S*q + q_old*diag(beta);
    q_old = q;
    for j=1:nrhs
      beta(j) = norm(q_new(:,j),2);
      if (beta(j) ~= 0) q(:,j) = q_new(:,j)/beta(j); end
    end
    theta = sqrt(theta1.*theta1 + beta.*beta);
    c_k = theta1./theta;
    s_k = beta./theta;
    delta = -s_old.*beta;
    theta1=  c_old.*theta;

    p = (q_old - p_old2*diag(delta_old))*diag(1./theta);

    x = x + p*diag(s)*diag(c_k);
    s = -s.*s_k;

    relres_hist = [relres_hist; abs(s)./r0];

    if (max(abs(s)./r0)<tol)
        its = i
        for j = 1:nrhs
            res(j) = norm(b(:,j)-alpha*x(:,j)-S*x(:,j), 2);
        end
        relres = res./r0;
        break;
    end

    p_old2 = p_old;
    p_old = p;
    s_old = s_k;
    c_old = c_k;
    delta_old = delta;

end
\end{PseudoCode}

\end{document}